\documentclass[10pt]{article}

\usepackage{amsfonts,latexsym,amstext}
\usepackage{amsmath}
\usepackage{amssymb}
\usepackage[english]{babel}
\usepackage[latin1]{inputenc}

\usepackage[usenames]{color}
\definecolor{red}{rgb}{1.0,0.0,0.0}

\definecolor{blu}{rgb}{0.0,0.0,1.0}

\definecolor{gre}{rgb}{0.03,0.50,0.03}

\definecolor{darkviolet}{rgb}{0.58, 0.0, 0.83}

\usepackage{enumerate}

\newtheorem{theorem}{Theorem}[section]
\newtheorem{lemma}[theorem]{Lemma}
\newtheorem{proposition}[theorem]{Proposition}
\newtheorem{definition}[theorem]{Definition}

\newtheorem{hypothesis}[theorem]{Hypothesis}

\newtheorem{remark}[theorem]{Remark}

\numberwithin{equation}{section}

\newcommand{\myref}[1]{(\ref {#1})}

\def\qed{{\hfill\hbox{\enspace${ \square}$}} \smallskip}
\def\sqr#1#2{{\vcenter{\vbox{\hrule height .#2pt \hbox{\vrule
 width .#2pt height#1pt \kern#1pt \vrule
width .#2pt} \hrule height .#2pt}}}}
\def\square{\mathchoice\sqr54\sqr54\sqr{4.1}3\sqr{3.5}3}

\def\qedo{\hbox{\hskip 6pt\vrule width6pt height7pt
depth1pt  \hskip1pt}\bigskip}

\def\eps{\varepsilon}

\def\ds{\begin{displaystyle}}
\def\eds{\end{displaystyle}}
\def\dis{\displaystyle }
\def\<{\left\langle }
\def\>{\right\rangle }

\def\dim{\noindent \hbox{{\bf Proof.} }}

\def\f{\mathbf f}

\def\R{\mathbb R}
\def\N{\mathbb N}

\def\E{\mathbb E}
\def\P{\mathbb P}

\def\cala{{\cal A}}

\def\cald{{\cal D}}

\def\calf{{\cal F}}

\def\calh{{\cal H}}

\def\call{{\cal L}}
\def\caln{{\cal N}}

\def\calu{{\cal U}}

\def\call{{\cal L}}

\def\calo{{\cal O}}

\def\Imm{{\operatorname{Im}}}

\def\to{\rightarrow}

\begin{document}

\title{Stochastic Control Problems with Unbounded Control Operators: solutions through generalized derivatives}
\date{}

 \author{Fausto Gozzi
\\
Dipartimento di Economia e Finanza,
Universit\`a LUISS - Guido Carli\\
Viale Romania 32,
00197 Roma,
Italy\\
e-mail: fgozzi@luiss.it\\
\\
Federica Masiero\\
Dipartimento di Matematica e Applicazioni, Universit\`a di Milano Bicocca\\
via Cozzi 55, 20125 Milano, Italy\\
e-mail: federica.masiero@unimib.it}

\maketitle
\begin{abstract}
This paper deals with a family of stochastic control problems
in Hilbert spaces which arises in many engineering/economic/financial applications (in particular the ones featuring boundary control and control of delay equations with delay in the control) and for which it is difficult to apply the dynamic programming approach due to the unboundedness of the control operator and to the lack of regularity of the underlying transition semigroup.
We introduce a specific concept of partial derivative, designed for this situation, and we develop a method to prove that the associated HJB equation has a solution with enough regularity to find optimal controls in feedback form.
\end{abstract}


\smallskip \noindent


\smallskip \noindent

\textbf{Acknowledgements}:
The authors thanks a lot the Associate Editor and the Referees for careful scrutiny and useful suggestions that led to an improved version of the paper.

Fausto Gozzi has been supported by the Italian
Ministry of University and Research (MIUR), in the framework of PRIN
projects 2015233N54 006 (Deterministic and stochastic evolution equations) and 2017FKHBA8 001 (The Time-Space Evolution of Economic Activities: Mathematical Models and Empirical Applications).

Federica Masiero has been supported by the Italian Ministry of University and Research (MIUR), in the framework of PRIN project 2015233N54 006 (Deterministic and stochastic evolution equations) and by the Gruppo Nazionale per l'Analisi Matematica, la Probabilit\`a e le loro Applicazioni (GNAMPA) of the Istituto Nazionale di Alta Matematica (INdAM).

\tableofcontents

\section{Introduction}

\subsection{Stochastic control in infinite dimension and its applications}

Stochastic optimal control problems arise in a large variety of applications (see e.g. \cite{Neck84}), and have been (and currently are) the object of an extensive theoretical and applied literature.
In recent years, due to progress in the methodologies and in the computational power there has been an increasing interest in studying also what is usually called ``the infinite dimensional case'', i.e. the case when the state and control variables take their values in infinite dimensional spaces. The infinite dimensional case allows for more realistic modeling and substantially means that the state equation is a Stochastic Partial Differential Equation (SPDE from now on) or a Stochastic Differential Delay Equation (SDDE from now on).

Such types of state equations arise naturally in a wide range of
applied models including physics, engineering, operations research, economics and finance.

On the one hand state equations of SPDE type are used when one wants to model control processes where the underlying dynamical system inherently depends on other basic variables beyond time.
For example we recall: the control of SPDEs arising in fluid dynamics (see e.g. \cite{DaPratoDebussche99}, \cite{Sritharanbook}),
in reaction-diffusion problems (see e.g. \cite{MouraFathy13}),
in modeling air pollution (see e.g. \cite{Do14AMS}, \cite{Seinfeld98}),
in robotics (see e.g. \cite{EvansPereiraBoutselisTheodoru19},
\cite{EvansKendallBoutselisTheodoru20}),
in elasticity theory and practice (see e.g. \cite{ChowMenaldi14} \cite{Do17Automatica},\cite{Do17JSV}),
in spatio-temporal economics growth models
(see e.g. \cite{BoucekkineCamachoFabbri13,Brito04}
in the deterministic case and
\cite{GozziLeocata21} in the stochastic case),
in advertising models (see e.g. \cite{BarucciGozziAOR99},
\cite{Huang12}).

On the other hand state equations of SDDE type are used when one wants to model control processes where the underlying dynamical system is not Markovian in the sense that the value of the state at time $t$ depends on the past of the state/control variables.
Such models arise in optimal advertising problems (see e.g. \cite{ChenWu20,GM,GMSJOTA,Huang12,Machowska19}); in optimal portfolio problems
(see e.g.  \cite{CarmonaEtAl18,FedFinSto,RosestolatoSwiech17});
in optimal production planning (see e.g.
 \cite{Chan07}, \cite{Yan00});
in the feedback stabilisation of engineering systems (see e.g. \cite{Li20AUTO}).

\subsection{The purpose of our paper: the ``unbounded control case'' and its interest}

Despite its interest for applications, the theory of stochastic control in infinite dimension is still a young and incomplete area. For this reason we think it is interesting to
continue to develop the theory for such problems, in particular
trying to cover new problems arising in applications like the ones quoted above.
For an up-to date account of the theory we send the reader see to the recent book \cite{FabbriGozziSwiech}; other books which partly look at the subject are
\cite[Chapter 13]{DP3}, \cite[Chapters 9-10]{Cerraibook}
and \cite[Chapters 5-6]{Nisiobook}).

This paper is devoted exactly to the above task in cases which are quite common in applications and for which there are, at the moment, few and incomplete results available: the cases, roughly speaking, where the control enters in the state equation in an ``unbounded'' way.
Such unboundedness translates into the fact that, while the state process takes values in a Hilbert space $H$
(e.g. $L^2(0,1)$), the term
which brings the control action into the state equation,
say $Cu$, takes values in a bigger Banach space
$\overline{H}$
(e.g. any Sobolev space of negative order) which strictly contains $H$.
Typical applications where such unboundedness arises are
\textbf{stochastic boundary control problems}
and \textbf{stochastic control of delay equations with delay in the control}.

Concerning \textbf{stochastic boundary control problems}, i.e. problems where, for structural reasons, the control is applied only at the boundary of a given region, are quite common in a wide range of applied problems.
In the engineering-related literature
we recall
some recent papers devoted to specific applications like, e.g., \cite{EvansKendallTheodoru21} in the field of robotics; \cite{Liu20AUTO,Lachemi20AUTO,Lachemi21TAC,MouraFathy13} for reaction-diffusion systems (also in the deterministic case);
\cite{Lamoline20} for stochastic port-Hamiltonian systems; \cite{Do17Automatica} \cite{Do17JSV} for marine risers and Timoshenko beams.
\\
Stochastic boundary control problems also arise in advertising models when one wants to take account of the age structure of the products, see e.g. \cite{BarucciGozziAOR99}, \cite{FaggianGrossetMMOR13}, \cite{GrossetViscolani},
\cite{Huang12}.

On the other hand, \textbf{stochastic optimal control problems with delay in the control}
arise since in many practical situations the effect of the control action persists in the future.
We recall, in this respect, the so-called carryover effect in advertising, see e.g. \cite{ChenWu20}, \cite{GM}, \cite{GMSJOTA}, \cite{Machowska19} and, for related deterministic problems, \cite{FeichtingerHartlSethi}.
\\
Applications to economics (delay in production due to the time to build) are studied in \cite{ChenWuAutomatica2010} where a problem with pointwise delay in the state and in the control is studied by means of the stochastic maximum principle (see also \cite{BambietalET} in a related deterministic case).
\\
Concerning applications to finance we recall \cite{CarmonaEtAl18} where a mean-field model of systemic risk with delay in the control is studied.
\\
Finally stochastic control problems with delay in the control are also related to the problem of information delay (i.e. the time which may be necessary
to implementing the control, which is studied e.g. in \cite{SZ}.

\subsection{The novelty of our methodology}

We use here the dynamic programming approach
with the aim of finding solutions of the associated HJB equations which are regular enough to write optimal control strategies in feedback form and to prove verification theorems.
In the literature one can, roughly speaking, distinguish three main methods to do this.
The first is to look at the theory of viscosity solutions (which is partly developed in these cases, see e.g.
\cite{FedericoGozziJDE,GozziRouySwiech06,Swiech20} and
\cite[Section 3.12 and 4.8.3]{FabbriGozziSwiech});
the viscosity solution is not differentiable in general, however in some infinite dimensional cases some type of differentiability can be proved
(see e.g. \cite{FedericoGoldysGozzi10SICON}, \cite{RosestolatoSwiech17}) however such methods seems not applicable here due to the unboundedness of the control operator.
The second approach (developed e.g. in \cite{ChowMenaldi97}, \cite{Ahmed01})
which we may call "variational" is based on the use of coercive bilinear forms associated to Ornstein-Uhlenbeck operators on a suitable Gelfand triple. Again this seems not applicable here since
it seems not compatible with the unboundedness of our control operator and the lack of null controllability estimates of our cases.

The third approach, which is it the one we use here is the theory of mild solutions (see e.g.
the papers \cite{CDP1}, \cite{CDP2}, \cite{DP3}, \cite{G1} \cite{Mas}, \cite{GoldysGozzi06SPA} and the book \cite[Chapter 4]{FabbriGozziSwiech}).
\\
Roughly speaking such approach, which works only in the semilinear case, rewrites the HJB equation
in a suitable integral form and tries to solve it using a fixed point argument. It is based on smoothing properties
of an underlying transition semigroup and allows to find solution which are regular enough to define optimal feedback control strategies.
\\
In the present case there are two technical barriers preventing the use of such approach:
\begin{itemize}
  \item
the need of giving a precise sense to the HJB equation, coming from the unboundedness explained above, which is the core of our setting;
  \item
the lack of smoothing properties of the underlying transition semigroup.
\end{itemize}
To deal with such issues we introduce
a specific concept of partial derivative, designed for this situation.
We observe that in the literature various concepts of infinite dimensional derivative has been used, depending on the context.
We mention, among others, the papers \cite{Gross1967},
\cite{LunardiRockner20} (which uses the so-called Fomin derivative), \cite{FTGgrad,Mas,Mas-inf-or} (which use the so-called $G$-gradient for a bounded operator $G$),
\cite[Chapter 4]{FabbriGozziSwiech}-\cite{FedericoGozziJDE}
(which use the so-called $G$-gradient for a possibly unbounded operator $G$).
Our definition extends the last one to take into account more
general cases of unboundedness, more precise explanations are given in Section \ref{subsection-C-directionalderivatives}.
\\
Once such derivative is introduced we perform a nontrivial extension of an idea (which we call ``partial smoothing'') that we used in our previous papers \cite{FGFM-I}-\cite{FGFM-II} in a delay case with bounded control operators. Details are given in Sections
\ref{sec:partsmooth-abstr-setting}-\ref{sec:convpartsmooth-abstr-setting}-\ref{sec-HJB}.

We must say that this paper is a first step to attack such difficult problems. Here we show how to find a regular solution to the HJB equation for a special type of cost functionals.
A second step, which is the object of our current research, is to cover more general cost functionals, in particular the ones where the current cost can be state dependent.
Moreover, to make our results more useful for applications we aim at proving prove verification type results: this would allow to construct the optimal controls in feedback form, on the line of what is done, e.g., in \cite{Ahmed03} or in our previous paper \cite{FGFM-II}.

We also think that it would be interesting
to study the case with partial observation
(on the line of what is done, in cases which do not include ours, in \cite{Ahmed15JMAA}, \cite{GozziSwiech00JFA}, \cite{BandiniCossoFuhrmanPham15AAP})
or the case of Mc Kean - Vlasov dynamics
(see e.g. \cite{Ahmed07}, or \cite{CossoEtAl20})
and to study the applicability, to our setting, of numerical methods for infinite dimensional HJB equations like the ones developed in \cite{AllaFalconeKalise15}.

\medskip

\subsection{Plan of the paper}

The plan of the paper is the following.
\begin{itemize}
  \item
  Section \ref{subsection-notation} introduces some basic notations.
   \item
   Section \ref{SE:EXCONTROL} introduces
   our driving examples showing how to rewrite them in
   a suitable infinite dimensional setting.
   \item
Section \ref{subsection-C-directionalderivatives}
provides the definition of our $C$-derivatives together with some comments to compare it with previous definitions.
   \item
   Section \ref{sec:partsmooth-abstr-setting} presents our partial smoothing result for Ornstein-Uhlenbeck semigroups
   (Proposition \ref{prop:partsmooth}).
   \item Section \ref{sec:convpartsmooth-abstr-setting}.
generalizes the partial smoothing result to convolutions
(Lemma \ref{lemma_convoluzione}).
\item
In Section \ref{sec-HJB} we first present a general control problem which includes both our driving examples (Subsection \ref{subsec-contr.pr-abstract}); then, in Subsection \ref{sub:HJBsol}, we state and prove our main result
of existence and uniqueness of mild solutions for the HJB equation (Theorem \ref{esistenzaHJB}).
\item
Appendix A is devoted to show that our motivating examples
satisfy the assumptions made.
\end{itemize}

\section{Basic notation and spaces}\label{subsection-notation}
For the reader's convenience we collect here the basic notation used throughout the paper.
\\
Let $H$ be a Hilbert space.
The norm of an element $x$ in $H$ will be denoted by
$\left|  x\right|_{H}$ or simply $\left|x\right|$,
if no confusion is possible, and by
$\left\langle \cdot,\cdot\right\rangle _H$,
or simply by $\left\langle \cdot,\cdot\right\rangle$ we denote the inner product in $H$.
We denote by $H^{\ast}$ the dual space of $H$.
If $K$ is another Hilbert space, $\call(H,K)$ denotes the
space of bounded linear operators from $H$ to $K$ endowed with the usual operator norm.
All Hilbert spaces are assumed to be real and separable.

Let $E$ be a Banach space. As for the Hilbert space case, the norm of an element $x$ in $E$ will be
denoted by $\left|x\right|_{E}$ or simply $\left|x\right|$,
if no confusion is possible.
We denote by $E^{\ast}$ the dual space of $E$,
and by $\left\langle \cdot,\cdot\right\rangle_{E^*,E}$
the duality between $E$ and $E^*$.

If $F$ is another Banach space, $\call(E,F)$ denotes the
space of bounded linear operators from $E$ to $F$ endowed with the usual operator norm.
All Banach spaces are assumed to be real and separable.

In what follows we will often meet inverses of operators which are not
one-to-one. Let $Q\in \call\left(H,K\right)$.
Then $H_{0}=\ker Q$ is a closed subspace of $H$. Let
$H_{1}:=[\ker Q]^\perp$ be the orthogonal complement of $H_{0}$ in $H$: $H_{1}$ is closed, too.
Denote by $Q_{1}$ the restriction of $Q$ to $H_{1}$: $Q_{1}$ is
one-to-one and $\operatorname{Im}Q_{1}=\operatorname{Im}Q$.
For $k\in \operatorname{Im}Q$, we define $Q^{-1}$ by setting
\[
Q^{-1}\left(k\right)  :=Q_{1}^{-1}\left(k\right)  .
\]
The operator $Q^{-1}:\operatorname{Im}Q\rightarrow H$ is called the pseudoinverse of $Q$. $Q^{-1}$ is linear and closed but in general not continuous.
Note that if $k\in\operatorname{Im}Q$, then
$Q_{1}^{-1}\left(  k\right)\in[\ker Q]^\perp$.
is the unique element of
\(
\left\{ h  :Q\left(  h\right)  =k\right\}
\)
with minimal norm (see e.g. \cite{Z}, p.209),

Next we introduce some spaces of functions.
Let $H$ and $Z$ be real separable Hilbert spaces.
By $B_b(H,Z)$ (respectively $C_b(H,Z)$, $UC_b(H,Z)$) we denote the space of all functions
$f:H\rightarrow Z$ which are Borel measurable and bounded (respectively continuous
and bounded, uniformly continuous and bounded).

Given an interval $I\subseteq \R$ we denote by
$C(I\times H,Z)$ (respectively $C_b(I\times H,Z)$)
the space of all functions $f:I \times H\rightarrow Z$
which are continuous (respectively continuous and bounded).
$C^{0,1}(I\times H,Z)$ is the space of functions
$ f\in C(I\times H, Z)$ such that, for all $t\in I$,
$f(t,\cdot)$ is continuously Fr\'echet differentiable with Fr\'echet derivative $\nabla f(t,x)\in \call(H,Z)$.
By $UC_{b}^{1,2}(I\times H,Z)$
we denote the linear space of the mappings $f:I\times H \to Z$
which are uniformly continuous and bounded
together with their first time derivative $f_t$ and their first and second space
derivatives $\nabla f,\nabla^2f$.
\\
If the destination space $Z=\R$ we do not write it in all the above spaces.
\\
The same definitions can be given if $H$ and $Z$ are Banach spaces.




\section{Two examples with unbounded control operator}
\label{SE:EXCONTROL}

We present here two stochastic controlled equations that motivates the introduction  of generalized partial derivatives in Section \ref{subsection-C-directionalderivatives}. What they have in common is that, once they are reformulated as infinite dimensional stochastic controlled evolution equations, the control operator is unbounded.

\subsection{Heat equations with boundary control}
\label{SSE:HEATEQUATION}

\subsubsection{The state equation}
\label{SSSE:SEBC}
Let $(\Omega, \calf,P)$ be a complete probability space endowed with a filtration $(\calf_t)_{t\geq 0}$. Fixed $0\le t \le T<+\infty$, we consider, in an open connected set with smooth boundary $\calo\subseteq \R^d$ ($d=1,2,3$) the controlled stochastic heat equation with Dirichlet boundary conditions and with boundary control:
\begin{equation}\label{eqDiri}
  \left\{
  \begin{array}{l}
  \dis
\frac{ \partial y}{\partial t}(s,\xi)
= \Delta y(s,\xi)+\dot{W}_Q(s,\xi), \qquad s\in [t,T],\;
\xi\in \calo,
\\\dis
y(t,\xi)=x(\xi),\; \xi\in \calo,
\\\dis
y(s,\xi)= u(s,\xi), \qquad s\in [0,T],\;
\xi\in \partial\calo.
\end{array}
\right.
\end{equation}
where $\Delta$ is the Laplace operator and we assume the following.

\begin{hypothesis}\label{hp:BC}
\begin{itemize}
\item[]
  \item[(i)]
The initial datum $x(\cdot)$ belongs to the state space $H:=L^2(\calo)$. The set $U$ of control values is a closed and bounded subset of the Hilbert space $K:=L^2(\partial\calo)$.
  \item[(ii)]
$\dot{W}_Q$ is a so-called colored space-time noise (with space covariance $Q\in \call(H)$), the filtration $(\calf_t)_{t\geq 0}$ coincides with the augmented filtration generated by $W_Q$;
  \item[(iii)]
the control strategy $u$ belongs to $\calu$ where
$$
\calu:=\left\lbrace u(\cdot):\Omega\times [0,T] \to U):\; \hbox{predictable}\right\rbrace
$$
\end{itemize}
\end{hypothesis}
Given any $(t,x) \in [0,T]\times \mathcal{O}$ and $u \in \mathcal{U}$, we denote, formally, by $y^{t,x,u}(s,\xi)$ the solution of \eqref{eqDiri} at $(s,\xi) \in [0,T]\times \mathcal{O}$.\footnote{Such solution could be defined with various methods. Here, similarly to
\cite[Appendix C]{FabbriGozziSwiech} we define such solution as the unique mild solution (see \eqref{eq:mildsolboundary} of the infinite dimensional system \eqref{eqDiri-abstr-contr} below.
Following the path outlined in
\cite[Appendix C]{FabbriGozziSwiech} we give sense to
\eqref{eqDiri} rewriting it as an evolution equation in the space $H:=L^2(\calo)$. We assume that the initial condition $x(\cdot)$ belongs to $H$.
The new state will be a process with values in $H$ given, formally, by $X(s;t,x,u)=y^{t,x,u}(s,\cdot)$.}
We define the operator $A_0$ in $H$ setting
(here $H^2(\calo)$ and $H^1_0(\calo)$ are the usual Sobolev spaces)
$$
\cald(A_0)=H^2(\calo)\cap H^1_0(\calo)
 \qquad
 A_0y = \Delta y
 {\rm \;\; for\; \;} y\in \cald(A_0).
 $$
The operator $A_0$ is self-adjoint and diagonal with strictly negative eigenvalues $\{-\lambda_n\}_{n\in \N}$
(recall that $\lambda_n\sim n^{2/d}$ as $n \to +\infty$).
We can endow $H$ with a complete orthonormal basis
$\{e_n\}_{n\in \N}$
of eigenvectors of $A_0$.\footnote{We know that $e_0$ is constant and, when $d=1$ and $\calo=(0,\pi)$, $(e_n(\xi))_{n\geq 1}:=(\sqrt{2}\sin (n\xi))_{n\geq 1}$.}
We recall that the linear trace operator
$D:L^2(\partial \calo)\rightarrow H$
is defined setting $D a=f$ where $f$ is the unique solution of the Dirichlet problem
$$
\left\{
\begin{array}{l}
\Delta f(\xi)=0,
\qquad \xi\in \calo,
\\
\dis
f(\xi)=a(\xi), \qquad \xi\in \partial\calo.
\end{array}
\right. $$
Equation (\ref{eqDiri}) can now be reformulated
(see \cite[Appendix C]{FabbriGozziSwiech} for a proof) as
\begin{equation}\label{eqDiri-abstr-contr}
  \left\{
  \begin{array}{l}
  \dis
d X(s)= A_0X(s)\,ds +
( -A_0)D u(s)\,dt+Q^{1/2}dW(s),
\\\dis
X(t)=x.
\end{array}
\right.
\end{equation}
where $W(\cdot)$ is a cylindrical noise in $H$.
We define
\begin{equation}\label{notazioneB}
B_0:=( -A_0)D.
\end{equation}
The operator $B_0$, defined in $K=L^2(\partial\calo)$,
does not take values in $H=L^2(\calo)$.
Indeed for all $\varepsilon>0$, the Dirichlet map takes its values in $\cald((-A_0)^{1/4-\varepsilon}$:
$D\in \call\left(K, \cald((-A_0)^{1/4-\varepsilon})\right)$
(see again \cite[Appendix C]{FabbriGozziSwiech}).
So,
$$
B_0=(-A_0)^{3/4+\varepsilon}(-A_0)^{1/4-\varepsilon}D:
K \rightarrow \cald((-A_0)^{-3/4-\varepsilon}):
$$
Here and from now on we take $0<\varepsilon<1/4$, indeed the point is to take $\varepsilon$ as small as possible in order to have in $B_0$ a better unbounded part.
Hence we have\footnote{Here, for $\gamma>0$ we denote by
$\cald((-A_0)^{-\gamma})$
the completion of $H$ with respect to the norm
$|\cdot|_{-\gamma}=|A_0^{-\gamma}\cdot|_H$.}
$B_0\in \call\left(K,\cald((-A_0)^{-3/4-\varepsilon})\right)$.
With an abuse of language with respect to the standard use, we
may say that $B$ is unbounded on $H$, in the sense that
its image is not contained in $H$ but in a space larger than $H$
(here
$\cald((-A_0)^{-3/4-\varepsilon}=H^{-3/2-2\varepsilon}(\calo)$)
which we will call $\overline{H}$.
\\
The unique mild solution (which exists thanks e.g. to
\cite[Theorem 1.141]{FabbriGozziSwiech}) of \eqref{eqDiri-abstr-contr}
is denoted by $X(\cdot;t,x,u)$ and is
\begin{equation}\label{eq:mildsolboundary}
X(s;t,x,u)=
e^{(s-t)A_0}x+\int_t^s e^{(s-r)A_0}B_0 u(r) dr +\int_t^se^{(s-r)A_0}Q^{1/2}dW(r)
,\text{ \ \ \ }s\in[t,T].
\end{equation}
Consequently, for any $(t,x) \in [0,T]\times \mathcal{O}$ and $u \in \mathcal{U}$, we give sense to
$y^{t,x,u}$ by setting $y^{t,x,u}(s,\xi):=X(s;t,x,u)(\xi)$.

\subsubsection{The optimal control problem}
\label{SSSE:OCBC}

For any given $t\in [0,T]$ and $x \in H$, the objective is to minimize, over all control strategies in $\calu$, the following finite horizon cost:
\begin{equation}\label{costoastratto}
J(t,x;u)=\E \left[\int_t^T \left[\ell_0(s)+\ell_1(u(s))\right]\,ds + \phi(X(T;t,x,u))\right],
\end{equation}
under the following assumption
\begin{hypothesis}\label{hp:BCcost}
\begin{itemize}
\item[]
  \item[(i)] $\ell_0:[0,T]\rightarrow \R$, is measurable and bounded.
  \item[(ii)] $\ell_1:U\to \R$ is measurable and bounded from below.
  \item[(iii)] $\phi:H\to \R$ is such that, for a suitable finite set $\{\f_1,\dots,\f_N\}\subseteq \cald((-A_0)^{\eta})$ (with $\eta>1/4$) and a suitable $\bar\phi \in B_b(\R^n)$ we have
      $$
      \phi_0(x)=\bar\phi\left(\<x,\f_1\>_H,\dots,
      \<x,\f_N\>_H\right).
      $$
\end{itemize}
\end{hypothesis}
Such cost can be seen as the rewriting in $H$ of a more ''concrete'' cost like
\begin{equation}\label{costoconcreto}
J_0(t,x;u)=\E \left[\int_t^T \int_{\mathcal{O}} \left[\bar\ell_0(s,\xi)+\bar\ell_1(u(s,\xi))\right]\,d\xi ds + \phi_0\left(y^{t,x,u}(T,\cdot)\right)\right],
\end{equation}
where $\bar\ell_0$ $\bar\ell_1$, $\phi_0$ are chosen so that the corresponding $\ell_0,\ell_1,\phi$ satisfy Hypothesis \ref{hp:BCcost} above.
\newline Note that here the current cost does not depend on the
state, this is due to the fact that putting the dependence on the state in the current cost would increase considerably the technical arguments in the solution of the HJB equation. Moreover , in order to bypass the lack of suitable smoothing properties of the underlying transition semigroup, we have to work on cost functionals which depend on the state only through a suitable operator $P$, which here turns out to be defined by projections on $\cald((-A_0)^{\eta})$ ,with $\eta>1/4$, see Section \ref{subsec-contr.pr-abstract} and \ref{sub:HJBsol} for further details.
\newline The value function of the problem is
\begin{equation}\label{valuefunction-diri}
 V(t,x):= \inf_{u \in \calu}J(t,x;u).
\end{equation}
We define the Hamiltonians leaving aside the term not depending on the control $u$.
For $p\in H$, $u \in U$, the current value Hamiltonian $\hat H_{CV}$ is given by (this is formal since neither $B_0u$ nor $B_0^*p$ belong to $H$, in general):
\begin{equation}
\label{eq:hatHCV}
\hat H_{CV}(p\,;u):=\<B_0u,p\>_{H}+\ell_1(u)=
\<u,B_0^*p\>_{H}+\ell_1(u)
\end{equation}
and the (minimum value) Hamiltonian by
\begin{equation}\label{psi1-diri}
\hat H_{min}(p):=\inf_{u\in U}\hat H_{CV}(p\,;u).
 \end{equation}
The associated HJB equation can then be formally written as
\begin{equation}\label{HJBformale-diri}
  \left\{\begin{array}{l}\dis
-\frac{\partial v(t,x)}{\partial t}=\cala [v(t,\cdot)](x) +\ell_0(t)+
\hat H_{min} (\nabla v(t,x)),\qquad t\in [0,T],\,
x\in \calh,\\
\\
\dis v(T,x)=\phi(x),
\end{array}\right.
\end{equation}
where $B$ is defined in {\eqref{notazioneB}}, and $\cala$ is the infinitesimal generator of the transition
semigroup $(R_{t})_{0 \leq t\leq T}$ associated to the process $X$ when the control is zero: namely $\cala$ is formally defined by
\begin{equation}\label{eq:ell-diri}
 \cala[f](x)=\frac{1}{2} Tr \; Q \nabla^2f(x)
+ \< x,A^*\nabla f(x)\>.
\end{equation}
From \eqref{eq:hatHCV} we easily see that, still formally,
$\hat H_{min}(p)$ dependsnot on $p $ but on $B_0^*p$. On the same line also the minimum point in \eqref{eq:hatHCV}, when it exists, only depends on $B_0^*p$. This means that the candidate optimal feedback map, if it exists, is a function of $B_0^*\nabla v$.
For this reason our main goal is to find a solution $v$ of \eqref{HJBformale-diri} for which $B_0^*\nabla v$ makes sense.
For this reason in the sequel we will use the notation
(for $p\in H$ and $q\in K$ such that the expressions below make sense):
\begin{equation}\label{eq:modham}
H_{CV}(q\,;u):=\<u,q\>_{K}+\ell_1(u)
\quad \hbox{and}\quad
H_{min}(q):= \inf_{u\in U} H_{CV}(q\,;u),
\end{equation}
so that
$$
\hat H_{CV}(p\,;u)=H_{CV}(B^*p\,;u)
\quad \hbox{and}\quad
\hat H_{min}(p)=H_{min}(B^*p)
$$

\subsection{SDEs with delay in the control}
\label{SSE:DELAYEQUATION}

\subsubsection{The state equation}
\label{SSSE:SEdelay}
In a complete probability space $(\Omega, \calf,  \P)$
we consider the following controlled stochastic
differential equation in $\R^n$ with delay in the control:
\begin{equation}
\left\{
\begin{array}
[c]{l}%
dy(s)  =a_0 y(s) ds+b_0 u(s) ds +\displaystyle\int_{-d}^0 u(s+\xi)b_1(d\xi) \, ds
+\sigma dW(s)
,\text{ \ \ \ }s\in [t,T] \\
y(t)  =y_0,\\
u(t+\xi)=u_0(\xi), \quad \xi \in [-d,0).
\end{array}
\right.  \label{eq-contr-rit}
\end{equation}
Here we consider the case of delay in the control, the case with delay also in the state is more complicated and cannot be treated as an application of the techniques introduced in the present paper.
\noindent We assume the following.

\begin{hypothesis}\label{hp:delaystate}
\begin{itemize}
\item[]
  \item[(i)] $W$ is a standard Brownian motion in $\R^k$, and $(\calf_t)_{t\geq 0}$ is the
augmented filtration generated by $W$;
  \item[(ii)] the control strategy $u$ belongs to $\calu$ where
$$\calu:=\left\lbrace u(\cdot):(\Omega\times [0,T]\to U):\;
\hbox{predictable} \right\rbrace
$$
with $U$ a closed and bounded subset of $\R^m$;
\item[(iii)] $a_0\in \call(\R^n;\R^n)$, $b_0 \in \call(\R^m;\R^n)$, $\sigma\in \call(\R^k;\R^n)$, $d>0$;
  \item[(iv)] $b_1$
  is an $m\times n$ matrix of signed Radon measures on $[-d,0]$ (i.e. it is an element of the dual space of $C([-d,0],\call(\R^m;\R^n))$).
\end{itemize}
\end{hypothesis}
Given any initial datum $(y_0,u_0)\in \R^n\times L^2([-d,0], \R^m)$ and any admissible control $u\in \calu$ equation (\ref{eq-contr-rit}) admits a unique strong (in the probabilistic sense) solution which is continuous and predictable
(see e.g. \cite{IkedaWatanabe} Chapter 4, Sections 2 and 3).

Notice that Hypothesis \ref{hp:delaystate}-(iv) on $b_1$ covers, but it is not limited to, the very common case of pointwise delay\footnote{Our delayed SDE includes, for example, the state equation used in \cite{CarmonaEtAl18}.} but it is technically complicated to deal with: indeed it gives rise, as we are going to see in next subsection, to an unbounded control operator $B$. We underline the fact that in the case of pintwise delay the matrix $b_1$ is a matrix of discrete measures, like weighted Dirac measures.
Recall that, in \cite{FGFM-I} and \cite{FGFM-II} the case of $b_1$ absolutely continuous with respect to the Lebesgue measure has been treated assuming
\begin{equation}\label{eq:b1restrictive}
b_1(d\xi)=\bar b_1(\xi)d(\xi),\;
\bar b_1\in L^2([-d,0],\call(\R^m;\R^n));
\end{equation}
it is clear that such an assumption on $b_1$ leaves aside the pointwise delay case which we treat here.

\subsubsection{Infinite dimensional reformulation}
\label{subsection-infdimref}

Now, using the approach of \cite{VK} (see \cite{GM} for the stochastic case), we reformulate equation (\ref{eq-contr-rit}) as an abstract stochastic differential equation in the Hilbert space $H=\R^n\times L^2([-d,0],\R^n)$.
To this end we introduce the operator $A_1 : \cald(A_1) \subset H \rightarrow H$ as follows: for $x=(x_0,x_1)\in H$, \begin{equation}\label{A1}
A_1x=( a_0 x_0 +x_1(0), -x_1'), \quad \cald(A_1)=\left\lbrace x\in H :x_1\in W^{1,2}([-d,0],\R^n), x_1(-d)=0 \right\rbrace.
\end{equation}
We denote by $A_1^*$ the adjoint operator of $A_1$:
\begin{equation}
 \label{Astar}
A_1^{*}x=( a_0^* x_0, x_1'), \quad \cald(A_1^{*})=\left\lbrace x\in H:x_1\in W^{1,2}([-d,0],\R^n), x_1(0)=x_0 \right\rbrace .
\end{equation}
We denote by $e^{tA_1}$ the $C_0$-semigroup generated by $A_1$. For $x\in H$ we have
\begin{equation}
e^{tA_1} \left(\begin{array}{l}x_0 \\x_1\end{array}\right)=
\left(
\begin{array}
[c]{ll}%
e^{ta_0 }x_0+\int_{-d}^{0}1_{[-t,0]} e^{(t+s)a_0 } x_1(s)ds \\[3mm]
x_1(\cdot-t)1_{[-d+t,0]}(\cdot).
\end{array}
\right)  \label{semigroup}
\end{equation}
Similarly, denoting by $e^{tA_1^*}=(e^{tA_1})^*$ the $C_0$-semigroup generated by $A_1^*$,
we have for
$z=\left(z_0,z_1\right)\in H$
\begin{equation}
e^{tA_1^*} \left(\begin{array}{l}z_0 \\z_1\end{array}\right)=
\left(
\begin{array}[c]{ll}
e^{t a_0^* }z_0 \\[3mm]
e^{(\cdot+t) a_0^* }z_0 1_{[-t,0]}(\cdot) +z_1(\cdot+t)1_{[-d,-t)}(\cdot).
\end{array}
\right)  \label{semigroupadjoint}
\end{equation}
The infinite dimensional noise operator is defined as
\begin{equation}
 \label{G}
G:\R^{k}\rightarrow H,\qquad Gy=(\sigma y, 0), \; y\in\R^k.
\end{equation}
The control operator $B_1$ is defined as
(here the control space is $K:=\R^m$ and we denote by $C'([-d,0],\R^n)$
the dual space of $C([-d,0],\R^n)$)
\begin{equation}
 \label{Bnotbdd}
\begin{array}{c}
B_1:\R^{m}\rightarrow \R^n \times C'([-d,0],\R^n),
   \\[2mm]
(B_1u)_0=b_0 u,
\quad \<f,(B_1u)_1\>_{C,C'} = \dis\int_{-d}^0\< f(\xi),b_1(d\xi )u\> , \quad u\in\R^m, \quad
f \in C([-d,0],\R^n).
\end{array}
\end{equation}
The adjoint $B_1^*$ is
\begin{equation}
 \label{B*notbdd}
\begin{array}{c}
B_1^*:\R^n \times C''([-d,0],\R^n) \rightarrow \R^{m},
   \\[2mm]
B_1^*(x_0,x_1)=
b^*_0 x_0+\dis\int_{-d}^0 b_1^*(d\xi)x_1(\xi)

,
\; (x_0,x_1)\in \R^n \times C([-d,0],\R^n),
\end{array}
\end{equation}
where we have denoted by $C''([-d,0],\R^n)$ the dual space of $C'([-d,0],\R^n)$, which contains $C([-d,0],\R^n)$: here we consider $B_1^*$ acting on $C([-d,0],\R^n)$, for a characterization of $C''([-d,0],\R^n)$, and for the inclusion of $C([-d,0],\R^n)$ in  $C''([-d,0],\R^n)$ see e.g. \cite{kaplan1}, \cite{kaplan2}, \cite{kaplan3} and \cite{Shannon}. If $b_1$ satisfies (\ref{eq:b1restrictive}), $B$ is a bounded operator from $\R^n$ to $H$, and we can easily write $e^{tA_1}B$: see \cite{FGFM-I}.
\\
If $b_1$ is as in Hypothesis \ref{hp:delaystate}-(iv), then $B$ is unbounded. Still it is possible to write $e^{tA_1}B$ by extending the semigroup, by extrapolation, to $\R^n\times C'([-d,0];\R^n)$.
We have, for $u \in \R^m$
\begin{equation}\label{eq:etAB}
\left(e^{tA_1}B_1\right)_0:\R^m \to \R^n,\qquad    \left(e^{tA_1}B_1\right)_0 u=
    e^{ta_0}b_0u+ \int_{-d}^0 1_{[-t,0]}e^{(t+r)a_0}b_1(dr)u,
\end{equation}
\begin{equation}\label{eq:etAB1}
\left(e^{tA_1}B_1\right)_1:\R^m \to C'([-d,0];\R^n),\qquad    \<f,\left(e^{tA_1}B_1\right)_1 u\>_{C,C'}=
    \int_{-d}^0 f(r+t) 1_{[-d,-t]}b_1(dr)u.
\end{equation}
Let us now define the predictable process
$Y=(Y_0,Y_1):\Omega \times [0,T]\to H$ as
$$
Y_0(s)=y(s), \qquad Y_1(s)(\xi)=\int_{-d}^\xi b_1(d\zeta)u(\zeta+s-\xi),
$$
where $y$ is the solution of \eqref{eq-contr-rit}
and $u\in \calu$ is the control process.
By \cite[Proposition 2]{GM}, the process $Y$
is the unique mild solution of the abstract evolution equation
in $H$
\begin{equation}
\left\{
\begin{array}
[c]{l}
dY(s)  =A_1Y(s) ds+B_1u(s) ds+GdW(s)
,\text{ \ \ \ }t\in[ 0,T] \\
Y(0)  =x=(x_0,x_1),
\end{array}
\right.   \label{eq-astr}%
\end{equation}
where $x_1(\xi)=\dis\int_{-d}^\xi u_0(\zeta-\xi)b_1(d\zeta)u_0(\zeta-\xi)$, for $\xi\in [-d,0)$, and $u_0$ has been introduced in \eqref{eq-contr-rit} as the initial condition of the control process. Since we have assumed $u_0\in L^2([-d,0], \R^m)$.
Note that we have $x_1\in L^2([-d,0];\R^n)$\footnote{This can be seen, e.g., by a simple application of Jensen inequality and Fubini theorem.}. The mild (or integral) form of (\ref{eq-astr}) is
\begin{equation}
Y(s)  =e^{(s-t)A_1}x+\int_t^se^{(s-r)A_1}B_1 u(r) dr +\int_t^se^{(s-r)A_1}GdW(r)
,\text{ \ \ \ }s\in[t,T]. \\
  \label{eq-astr-mild}%
\end{equation}
Here, similarly to what happen in the previous example (see Subsection \ref{SSSE:SEBC}), we may say that the image of $B$ is not contained in $H$ but in a space larger than $H$ (here $\R^n\times C'([-d,0],\R^n)$) which we will call $\overline{H}$.

\subsubsection{The optimal control problem}
\label{SSSE:OCdelay}

Similarly to the previous section the objective is to minimize, over all control strategies in $\calu$, a finite horizon cost:
\begin{equation}\label{eq:costoconcretodelay}
\bar J(t,y_0,u_0;u(\cdot))=\E \left[\int_t^T \left[\ell_0(s)+\ell_1(u(s))\right]\,ds + \bar\phi(y(T;t,x))\right]
\end{equation}
undert the following assumption
\begin{hypothesis}\label{hp:delaycost}
\begin{itemize}
\item[]
  \item[(i)] $\ell_0:[0,T]\rightarrow \R$, is measurable.
  \item[(ii)] $\ell_1:U\to \R$ is measurable and bounded from below
  \item[(iii)] $\bar\phi:\R^n\to \R$ is measurable and bounded.
\end{itemize}
\end{hypothesis}
Such cost functional, using the infinite dimensional reformulation given above, can be rewritten as
\begin{equation}\label{eq:costoastrattodelay}
J(t,x;u(\cdot))=\E \left[\int_t^T \left[\ell_0(s)+\ell_1(u(s))\right]\,ds +
\phi(Y(T;t,x))\right]
\end{equation}
where $\phi:H\to \R$ is defined as $\phi(x_0,x_1)=\bar\phi(x_0)$
for all $x=(x_0,x_1)\in H$. Note again that here the current cost does not depend on the
state, again this is due to the fact that putting the dependence on the state in the current cost would increase considerably the technical arguments in the solution of the HJB equation.
\newline The value function of the problem is
\begin{equation}\label{valuefunction-delay}
 V(t,x):= \inf_{u \in \calu}J(t,x;u).
\end{equation}
The Hamiltonians can be defined exactly in the same way as in Subsubsection \ref{SSSE:OCBC} and (using the modified Hamiltonians introduced in \eqref{eq:modham})
the associated HJB equation is formally written as
\begin{equation}\label{HJBformale-delay}
  \left\{\begin{array}{l}\dis
-\frac{\partial v(t,x)}{\partial t}=\cala [v(t,\cdot)](x) +\ell_0(t)+
H_{min} (B_1^*\nabla v(t,x)),\qquad t\in [0,T],\,
x\in H,\\
\\
\dis v(T,x)=\phi(x),
\end{array}\right.
\end{equation}
where $B_1$ is defined in {\eqref{Bnotbdd}}, and $\cala$ is the infinitesimal generator of the transition
semigroup $(R_{t})_{0 \leq t\leq T}$ associated to the process $Y$ when the control is zero: namely $\cala$ is formally defined by
\begin{equation}\label{eq:ell}
 \cala[f](x)=\frac{1}{2} Tr \; GG^* \nabla^2f(x)
+ \< x,A_1^*\nabla f(x)\>.
\end{equation}
On the same line of Subsubsection \ref{SSSE:OCBC} the candidate optimal feedback map, if it exists, is a function of $B_1^*\nabla v$.


\section{$C$-derivatives}\label{subsection-C-directionalderivatives}

In this Section we introduce the definition of generalized partial derivatives
(that we call $C$-directional derivatives, where $C$ is a suitable linear operator)
which is suitable for our needs.
$C$-directional derivatives of functions have been introduced in
\cite[Section 2]{Mas}, \cite{FTGgrad} in the case when $C$ is a bounded operator (see also \cite{FGFM-I,FGFM-II}),
and in \cite{FedericoGozziJDE},
\cite[Section 4.2.1]{FabbriGozziSwiech} in the case when $C$ is possibly unbounded.

Our definition is different from the ones recalled above and is designed to cover a wider class of ``unbounded'' examples, in particular it makes it possible to treat the case when the intersection of image of $C$ and the state space is just the origin, which is, e.g., the case of our examples of Section \ref{SE:EXCONTROL} which were not treatable within the previous setting.

We also recall that concepts which are connected to the one of $C$-directional derivative are the one of Fomin derivative
(see, e.g., \cite[Chapter 3]{BogachevAMS2010} and, recently, \cite{LunardiRockner20}) and the one of derivative
in the directions of a proper subspace
(see e.g. \cite{Gross1967}).

Here is our new definition. The operator $C$ is still ``possibly unbounded'' in the sense that it does not take its values in the state space $H$ but in a larger Banach space $\overline H$ such that $H\subset \overline H$ with continuous embedding.

\begin{definition}
\label{df4:Gderunbounded} Let $H, \,Z,\, K$ and $\overline H$ be
Banach spaces such that $H\subset \overline H$ with continuous embedding.
Let $C:K\rightarrow \overline H$ be a linear and bounded operator.
\begin{itemize}
\item[(i)]
Let $k\in K$ and let $f:\overline H\rightarrow Z$.
We say that $f$ admits $C$-directional derivative
at a point $x\in \overline{H}$ in the direction $k\in K$
(and we denote it by $\nabla^{C}f(x;k)$) if the limit, in the norm topology of $Z$,
\begin{equation}\label{Cderivatabis}
 \nabla^{C}f(x;k):=\lim_{s\rightarrow 0}
\frac{f(x+s Ck)-f(x)}{s},
\end{equation}
exists.
\item[(ii)]
Let $f:\overline H\rightarrow Z$.
We say that $f$ is $C$-G\^ateaux differentiable
at a point $x\in \overline{H}$ if $f$ admits the $C$-directional derivative in every
direction $k\in K$ and there exists a {\bf bounded} linear operator,
the $C$-G\^ateaux derivative $\nabla^C f(x)\in \call(K,Z)$, such that $\nabla^{C}f(x;k)  =\nabla^{C}f(x)k$
for all $k \in K$. We say that $f$ is $C$-G\^ateaux
differentiable on $H$ (respectively $\overline{H}$) if it is $C$-G\^ateaux differentiable at every point $x\in H$
(respectively $x\in\overline{H}$).
\item[(iii)]
Let $f:\overline H\rightarrow Z$.
We say that $f$ is $C$-Fr\'echet differentiable
at a point $x\in \overline{H}$ if it is $C$-G\^ateaux differentiable and if the limit
in (\ref{Cderivatabis}) is uniform for $k$ in the unit ball of $K$.
In this case
we call $\nabla^C f(x)$ the $C$-Fr\'echet derivative (or simply the $C$-derivative) of $f$ at $x$. We say that $f$ is $C$-Fr\'echet differentiable on $H$ (respectively $\overline{H}$) if it is $C$-Fr\'echet differentiable at every point $x\in H$ (respectively $x\in\overline{H}$).
\end{itemize}
\end{definition}

\medskip
 \begin{remark}
 \label{rm:Gderunbounded1}
{\rm
The main idea behind the use of $C$-derivatives (starting from the papers \cite{Mas} and \cite{FTGgrad}) lies in the fact that, in applying the dynamic programming approach to optimal control problems which are linear in the control (with control operator $C:U \to H$ where $U$ is the control space and $H$ is the state space), the natural regularity requirement needed on the value function $V$ to write the optimal feedbacks is that $\nabla^C V$ is well defined.\footnote{Here we are simplifying a bit since, as one can read in
\cite[Section 4.8.1.4]{FabbriGozziSwiech} (in particular equation (4.294)), the operator $C$ in the gradient may be chosen a bit differently, and the linearity in the control can be weakened without affecting the main issues.} This means that only directional derivatives in the directions of the image of $C$ matter for the purpose of writing optimal feedback controls. In some cases, like the distributed control of heat equation (see e.g. \cite[Section 2.6.1]{FabbriGozziSwiech}), the image of $C$ is contained in the state space (call it $H$), so $\nabla^C V$ is always well defined when $\nabla V$ exists. In some other cases, like the boundary control or the pointwise delayed control (see e.g. \cite[Sections 2.6.2 and 2.6.8]{FabbriGozziSwiech}) the image of the control operator $C$ is not contained in $H$ and it may even happen that the intersection of this image with $H$ is only the origin, which is the  case of the driving examples of this paper.

One strategy, used e.g. in \cite{FedericoGozziJDE,FedericoGozziAAP} and in \cite[Section 4.8]{FabbriGozziSwiech} to deal with such cases is to decompose the control operator $C$ in the product $C_1C_2$ where $C_2:K\to H$ is bounded while the ``unbounded part'' $C_1$ is a closed unbounded operator $C_1:D(C_1)\subseteq H  \to H$ which usually is a power of the operator $A$ driving the state equation.
In this case the derivative needed to express the feedback
control is $\nabla^{C_1} V$ which is defined exactly as
in \cite[Definition 2.2]{FedericoGozziJDE} or
\cite[Definition 4.4]{FabbriGozziSwiech}.\footnote{Note that \cite[Definition 4.4]{FabbriGozziSwiech} is more general
than our Definition in the sense that it allows the
operator $C$ to depend on the state variable $x \in H$. This could be performed here with ideas similar to what is done in \cite[Section 4.2]{FabbriGozziSwiech}.
We do not do this since it would increase the technicalities
without changing the main ideas which we want to make clear for the reader.}

In such setting, due to the boundedness required
in \cite[Definition 2.2-(ii)]{FedericoGozziJDE}, asking that $\nabla^{C_1} V$ exists
substantially means that we consider the directional derivatives
of $V$ in the directions of $\operatorname{Im} \bar C_1$ where $\bar C_1$ is the
extension of $C_1$ from the whole $H$ to a suitable extrapolation space.
The image of $\bar C_1$ contains
(but can be much larger than) the one of the control
operator $C$\footnote{For example, in the case of Neumann or Dirichlet boundary control in dimension 1, the image of $C$ is two-dimensional while the one of $\bar C_1$ is infinite dimensional}.

The approach used here is sharper in the sense that we look exactly at the derivatives in the directions of the image of $C$, even if they go out of the state space $H$. In this way we also avoid working with the decomposition of the operator $C$, which is not sharp for our purposes, in particular in the case of pointwise delayed control of Subsection \ref{SSE:DELAYEQUATION} since in this case fractional powers of $A$ are not well defined.\footnote{A similar issue would arise if we consider boundary control problems where the driving operator $A$ is of first order, like in the case of age-structured problems, see e.g., in the deterministic case, \cite{FaggianGozziKort}.}
}
\hfill\qedo
\end{remark}

 \begin{remark}
 \label{rm:Gderunbounded2}
\rm
Definition \ref{df4:Gderunbounded} is exactly the definition of
$C$-derivative contained in \cite[Section 2]{Mas}, \cite{FTGgrad} when
$\overline{H}=H$. This means
that, in this case, the classical G\^ateaux or Fr\'echet differentiability
implies the $C$-G\^ateaux or $C$-Fr\'echet differentiability.

In \cite[Definition 4.4]{FabbriGozziSwiech} the operator $C$ is a closed, possibly unbounded, linear operator
$C:D(C)\subseteq K \to H$.
This case can be partly embedded in the one we consider in Definition \ref{df4:Gderunbounded}.
We explain now why, restricting to the case when $H$ is reflexive, which is true in our examples.

Let $C^*:D(C^*)\subseteq H' \to K'$ be the adjoint of $C$ defined in the usual way through the duality $\<C^*h,k\>_{K',K}=\<h,Ck\>_{H',H}$, $\forall\,k\in D(C) ,\, \forall\,h\in D(C^*) $.
By \cite[Theorem 5.29]{Kato76}, since $H$ is reflexive, we know that $C^*$ is densely defined.
Let
\begin{equation}\label{eqnuova}
E:=D(C^*)=\left\lbrace e \in H' : \exists \,  a>0:\,\forall\, k\in D(C)\, \vert\< Ck, e\>_{H,H'}\vert\leq a|e|_{H'}\right\rbrace\subseteq H',\end{equation}
 endowed with the usual graph norm, i.e.
$$
\|w\|_E:=\|w\|_{H'}+\|C^*w\|_{K'}, \qquad \forall w \in E.
$$
Let then $E':=D(C^*)'$. Clearly by \eqref{eqnuova} duality
$H''\subseteq E'$.
Then, by the canonical embedding of the bidual we have
$H \subseteq H''\subseteq E'$.
We extend, by extrapolation (see e.g.
\cite[\S II.5]{EngelNagelBook} for the general theory and
\cite[\S 3.3]{FaggianGozziKort} or \cite{Faggian2005,FaggianDCDIS} for specific cases)
$C$ to a continuous operator $\widetilde C:K\to E'$ setting, for $k \in K$ and $y\in E$,
$$
\<\widetilde Ck,y\>_{E',E}= \<k,C^*y\>_{K,K'}
$$
Continuity of $\widetilde C$ immediately follows observing that
$$
|\<\widetilde Ck,y\>_{E',E}|\le |\<k,C^*y\>_{K,K'}|
\le |k|_K |C^*y|_{K'}
$$
and taking the supremum over all $y\in E$ in the unit ball.\footnote{
Notice that the second adjoint operator $C^{**}: D(C^{**})\subset K''\rightarrow H''$ is defined through the equality:
$$
 \<k,C^*y\>_{K,K'}= \<C^{**}k,y\>_{H'',H'}\;\forall\,k\in D(C^{**}) ,\, \forall\,y\in D(C^*) ;
$$
with $D(C^{**})$ defined analogously to \eqref{eqnuova}. So $ \widetilde C$ and $C^{**}$ are operator acting and taking values on different spaces:
$$
\<\widetilde Ck,y\>_{E',E}=\<C^{**}k,y\>_{H'',H'}.
$$
}
\newline In this context we now compare
\cite[Definition 4.4]{FabbriGozziSwiech}
for $C$ and Definition \ref{df4:Gderunbounded} for the corresponding extension $\widetilde C$.
Indeed we observe that
\cite[Definition 4.4-(i)]{FabbriGozziSwiech}
says, at point (i) (definition of directional derivatives):

\emph{``The $C$-directional
 derivative of $f$ at a point $x\in H$ in the direction $k\in D(C)\subseteq K$ is defined as:
 \begin{equation}
 \nabla^{C}f(x;k):=\lim_{s\rightarrow 0}
\frac{f(x+s Ck)-f(x)}{s},\text{ }s\in\mathbb{R},
 \label{Cderivata}
 \end{equation}
provided that the limit exists.''}

But $k\in D(C)$, in the setting introduced above
means that $\widetilde Ck\in H$.
Hence, concerning point (i), Definition
\ref{df4:Gderunbounded} extends
\cite[Definition 4.4]{FabbriGozziSwiech}.

Finally we observe that, when the image of the operator $C$ crosses
$H$ only at the origin, then,
\cite[Definition 4.4]{FabbriGozziSwiech} cannot
be used while Definition \ref{df4:Gderunbounded} is still fit.
\hfill\qedo
\end{remark}

\medskip
 \begin{remark}
 \label{rm:Gderunbounded3}
{\rm
Observe that, similarly to what observed in \cite[Remark 4.5]{FabbriGozziSwiech} for Definition 4.4 (see also \cite[Definition 2.2]{FedericoGozziJDE}), even if $f$ is Fr\'echet differentiable at $x \in H$, the $C$-derivative may not exist in such point.
This is obvious if we take, e.g., $f(x)=|x|^2$, $C:K\to \bar H$
with $\operatorname{Im} C \not\subseteq H$. If $k\in K$ is such that $Ck\not\in H$
clearly $\nabla^{C}f\left( x;k\right)$ does not exist.}
\hfill\qedo
\end{remark}

We are now in position to define suitable spaces of $C$-differentiable functions.

\begin{definition}
\label{df4:Gspaces}
Let $I$ be an interval in $\R$, let $H$, $\overline{H}$, $K$ and $Z$ be suitable real Banach spaces. Moreover let $H\subset \overline{H}$ with continuous inclusion, and let $C\in \call(K,\overline{H})$.
\begin{itemize}
\item
We call $C^{1,C}_{b}(\overline{H},Z)$ the space of all continuous and bounded functions $f:\overline H\to Z$ which admit continuous and bounded $C$-Fr\'echet derivative.
Moreover we call $C^{0,1,C}_b(I\times \overline{H},Z)$ the space of
continuous and bounded functions $f:I\times \overline H\to Z$ such that
for every $t\in I$, $f(t,\cdot)\in C^{1,C}_b(\overline{H},Z)$ and
$\nabla^C f\in C_b\left(I\times \overline{H},L(K,Z)\right)$.
When $Z=\R$ we write  $C^{1,C}_{b}(\overline{H})$ instead of $C^{1,C}_{b}(\overline{H},Z)$, and it turns out that if $f\in C^{1,C}_{b}(\overline{H})$, then  $\nabla^C f\in C_b\left(I\times \overline{H},K')\right)$.
\item
For any $\alpha\in(0,1)$ and $T>0$ (this time $I$ is equal to $[0,T]$) we denote by
$C^{0,1,C}_{\alpha}([0,T]\times \overline{H},Z)$
the space of functions
$f\in C_b([0,T]\times H,Z)\cap
C^{0,1,C}_b((0,T]\times \overline{H},Z)$\footnote{Note that here
$f(t,\cdot)$ is well defined only in $H$ when $t=0$, while for $t>0$ it is defined over $\overline{H}$. The reason is that the Ornstein-Uhlenbeck semigroup in our examples and in our setting (and consequently the solution of the HJB equation) satisfy the same property.}
such that
the map $(t,x)\mapsto t^{\alpha} \nabla^C f(t,x)$
belongs to $C_b((0,T]\times \overline{H},\call(K,Z))$.
When $Z=\R$ we omit it.
The space $C^{0,1,C}_{\alpha}([0,T]\times \overline{H},Z)$
is a Banach space when endowed with the norm
\[
 \left\Vert f\right\Vert _{C^{0,1,C}_{\alpha}([0,T]\times \overline{H},Z)  }=\sup_{(t,x)\in (0,T]\times \overline{H}}
\vert f(t,x)\vert+
\sup_{(t,x)\in (0,T]\times \overline{H}}  t^{\alpha }\left\Vert \nabla^C f(t,x)\right\Vert_{\call(K,Z)}.
\]
When clear from the context we will write simply
$\left\Vert f\right\Vert _{C^{0,1,C}_{\alpha}}$.
\end{itemize}
\end{definition}


%


\section{Partial smoothing for Ornstein-Uhlenbeck semigroups}
\label{sec:partsmooth-abstr-setting}

In this section we study the ``partial smoothing'' properties of the Ornstein-Uhlenbeck semigroup (which we call $R_t$, for $t\ge 0$) applied to a generic function $f$ weakening the definition of ``smoothing'' given, e.g., in \cite{DPZ91} (see also \cite[Chapter 9]{DP1}).
Note that a type of partial smoothing has been already developed, e.g., in \cite[Ch.4]{FabbriGozziSwiech} and in \cite{FGFM-I,FGFM-II}.
As said above, the main difference here is that the directions along which we take the derivative can go out of the state space $H$ and this allows to treat in sharper way the control problems exposed in Section \ref{SE:EXCONTROL}.

The following basic assumption holds throughout this section.
\begin{hypothesis}\label{ip-sde-common}
\begin{enumerate}[(i)]
\item[]
\item
Let $H$, $K$ $\Xi$ be three real separable Hilbert
spaces\footnote{These will be usually the state space,
the control space and the noise space, respectively.}.
  \item
Let $(\Omega, \calf,(\calf_t)_{t\geq 0}, \P)$ be a filtered probability space
satisfying the usual conditions and let $W$ be an
$(\Omega, \calf,(\calf_t)_{t\geq 0}, \P)$-cylindrical Wiener process in $\Xi$ where $(\calf_t)_{t\geq 0}$ is the augmented filtration generated by $W$.
  \item
Let $A:D(A)\subseteq H \to H$ be the generator of a strongly continuous semigroup
$e^{tA},\, t\geq 0$ in $H$,
  \item
Let $G\in\call(\Xi,H)$ be such that the selfadjoint operator
\begin{equation}\label{cov-gen}
Q_t=\int_0^t e^{sA}GG^*e^{sA^*}\,ds
\end{equation}
is trace class. We call $Q=GG^*\in \call(H)$.
\end{enumerate}
\end{hypothesis}

Let $Z(\cdot;x)$ be the Ornstein-Uhlenbeck process
which solves the following SDE in $H$.
\begin{equation}\label{ornstein-gen}
\left\lbrace\begin{array}{l}
dZ(t)=AZ(t)dt+GdW(t),\\
X(0)=x.
\end{array}\right.
\end{equation}
The process $Z(\cdot;x)$ is to be considered in its mild formulation:
\begin{equation}
Z(t;x)  =e^{tA}x +\int_0^te^{(t-s)A}GdW(s)
,\text{ \ \ \ }t\ge 0. \\
  \label{ornstein-mild-gen}
\end{equation}
$Z$ is a Gaussian process, namely for every $t>0$, the law of
$Z(t)$ is $\caln (e^{tA}x,Q_t)$, the Gaussian measure with mean $e^{tA}x$ and
covariance operator $Q_t$ defined in (\ref{cov-gen}).
The convolution $\int_0^te^{(t-s)A}GdW_s$ has law
$\caln (0,Q_t)$ and will be sometimes denoted by $W_A(t)$.
The associated Ornstein-Uhlenbeck transition semigroup $R_t$ is defined by setting, for every $\psi\in B_b(H)$ and $x\in H$,
\begin{equation}
 \label{ornstein-sem-gen}
R_t[\psi](x)=\E \psi(Z(t;x))
=\int_H \psi(z+e^{tA}x)\caln(0,Q_t)(dz).
\end{equation}
To study regularizing properties in the directions of an ``unbounded'' operator $C$ (as introduced in Section \ref{subsection-C-directionalderivatives}), and for functions that have a special dependence on the state, through an operator $P$ that we are going to introduce, we assume the following.

\begin{hypothesis}\label{ip:PC}
\begin{enumerate}[(i)]
\item[]
\item
Let $\overline H$ be a real Banach space such that
$H\subseteq \overline{H}$ with continuous and dense inclusion
and that the semigroup $e^{tA}$ admits an extension $\overline{e^{tA}}:\overline{H}\to \overline{H}$ which is still a $C_0$ semigroup.
\item
Let $C\in \call(K, \overline H)$.
\item Let $P:H \to H$ be linear and continuous.
Assume that, for every $t>0$ the operator $Pe^{tA}:H \to H$ can be extended to a continuous linear operator $\overline H\to H$, which will be denoted by $\overline{Pe^{tA}}$. With this notation the operator $\overline{Pe^{tA}}C:K\to H$ is well defined and continuous.
\end{enumerate}
\end{hypothesis}

We now provide two remarks on the above hypothesis: the first on the adjoint of $Pe^{tA}$, the second one on the validity of such hypothesis in our examples.

\begin{remark}\label{rm:adjointPeta}
In the framework of the above Hypothesis \ref{ip:PC}
it is natural to identify $H$ with its topological dual $H'$
and consider the Gelfand triple
$$
\overline{H}'\subseteq H\subseteq\overline{H}.
$$
The adjoint
of the operator $Pe^{tA}:H\to H$ which is, clearly,
$e^{tA^*}P^*:H\to H$, indeed takes its values in
$\overline{H}'$ and is, consequently, the adjoint
$$
\left(\overline{Pe^{tA}}\right)^*: H\to \overline{H}'
$$
of the extended operator $\overline{Pe^{tA}}$.
\\
Indeed, consider $\{x_n\}\subset H$ such that, in the topology of $\overline{H}$, we have $x_n \to \bar x\in \overline H$.
We know, by Hypothesis \ref{ip:PC}, that
$Pe^{tA}x_n \to \overline{Pe^{tA}}\bar x$, hence,
for every $y \in H$,
$$
\<x_n,e^{tA^*}P^*y\>_H=\<Pe^{tA}x_n,y\>_H \to \<\overline{Pe^{tA}}\bar x,y\>_H.
$$
Hence, the continuous linear form $\pi_y$ on $H$ (represented, with the Riesz identification on $H$, by $e^{tA^*}P^*y$) given by
$$
\pi_y:H\to \R, \qquad \pi_y (h)=\<h,e^{tA^*}P^*y\>_H, \qquad h\in H,
$$
can be extended to a continuous linear form
$$
\bar \pi_y:\overline{H} \to \R, \qquad
\bar\pi_y (\bar h)=\<\overline{Pe^{tA}}\bar h,y\>_H
\qquad \bar h\in \overline{H},
$$
with $|\bar\pi_y (\bar h)|\le |\overline{Pe^{tA}}|_{\call(\overline{H},H)}
|\bar h|_{\overline{H}} |y|_H$.
This is equivalent to say that, under the Riesz identification of $H$ with $H'$,
$e^{tA^*}P^*y\in \overline{H}'\subseteq H$.
\end{remark}

\begin{remark}\label{rm:computeexamples}
In the case of Subsection \ref{SSE:HEATEQUATION}
the above Hypotheses \ref{ip-sde-common} and \ref{ip:PC}
are satisfied if we choose, as seen in Subsection
\ref{SSE:HEATEQUATION},
$$
H=L^2(\calo), \qquad
\overline H= \cald\left((-A_0)^{-3/4-\eps}\right)
=H^{-3/2-2\eps}(\calo)
\quad \hbox{(for suitable small $\eps>0$),}
$$
$A=A_0$, $C=B=(-A_0)D$ as from \eqref{notazioneB}, $P$ any continuous operator $H\to H$
(we will later take $P$ to be a finite dimensional projection whose image is contained in
$\cald\left((-A_0)^{-\eta}\right)$ for some $\eta\ge 0$).
Since we can extend immediately $e^{tA_0}$ to
$$
\overline{e^{tA_0}}:\overline H \to  H
$$
then, in this case, $\overline{Pe^{tA}}=P\overline{e^{tA}}$.

In the case of Subsection \ref{SSE:DELAYEQUATION}
the above Hypotheses \ref{ip-sde-common} and \ref{ip:PC}
are satisfied if we choose
$$
H=\R^n \times L^2(-d,0;\R^n),\qquad
\overline H=\R^n \times C'([-d,0];\R^n),
$$
(but also $\overline H=\R^n \times W^{-1,2}([-d,0];\R^n)$ can be chosen),
$A=A_1$, $C=B$ as from \eqref{Bnotbdd}, and $P(x_0,x_1)=(x_0,0)$. Here the embedding of $L^2([-d,0];\R^n)\subset C'([-d,0];\R^n)$ is to be considered in the following sense: to any $f\in L^2([-d,0];\R^n)$ we associate the measure $\mu_f\in C'([-d,0];\R^n)$ such that $\mu_f(d\xi)=f(\xi)d\xi$.
\newline Note that
$$
\operatorname{Im }P=\R^n \times \{0\}.
$$
Moreover, by \myref{semigroup}, we have, for $x=(x_0,x_1) \in H$,
$$
Pe^{tA}x=\left(
e^{ta_0 }x_0+\int_{-d}^{0}1_{[-t,0]} e^{(t+s)a_0 } x_1(s)ds,
0 \right)
$$
Hence, also in this case, we can extend immediately $Pe^{tA}$ to
$$
\overline{Pe^{tA}}:\overline H \to H
$$
by setting, for $x=(x_0,x_1) \in \R^n \times C'([-d,0];\R^n)$
\begin{equation}
\overline{Pe^{tA}}x=\left(
e^{ta_0 }x_0+\int_{-d}^{0}1_{[-t,0]} e^{(t+s)a_0 } x_1(ds),
0 \right).
\label{PetAbardelay}
\end{equation}
Hence, also here Hypothesis \ref{ip:PC}-(iii) is satisfied.

Notice that in the second example $P$ can be immediately extended to $\overline{P}:\overline{H}\to H$ so $\overline{Pe^{tA}}=\overline{P}\,\,\overline{e^{tA}}$
while in the first example $P$ may not admit such an extension
(it does when $P$ is a finite dimensional projection).

Finally notice that in both examples we have
$\Imm  \overline{Pe^{tA}}\subseteq \Imm P$.
\hfill\qedo
\end{remark}



We pass to define the spaces where our ``initial'' data will belong.

\begin{definition}\label{df:spaziphi1}
We call $B_b^P(H)$ (respectively $C_b^P(H)$, $UC_b^P(H)$) the set of functions $\phi:H\to \R$ for which there exists
$\bar\phi : \operatorname{Im}(P)\to \R$
bounded and Borel measurable and (respectively continuous, uniformly continuous)\footnote{Here we endow $Im P\subseteq H$ with the topology inherited by $H$.} such that
\begin{equation}\label{fi-gen-allargata}
\phi(x)=\bar\phi(Px) \quad
\forall x\in H.
\end{equation}
\end{definition}


\begin{remark}\label{rm:ipPCunifcont}
We observe that, in the above Definition \ref{df:spaziphi1}, when $\bar\phi : \operatorname{Im}(P)\to \R$
is Borel measurable (respectively continuous, uniformly continuous), then also $\phi$ is Borel measurable (respectively continuous, uniformly continuous). Hence we can easily see that $B^P_b(H)$ (respectively $C^P_b(H)$, $UC^P_b(H)$) is a linear subspace of $B_b(H)$ (respectively $C_b(H)$, $UC_b(H)$).
%

We also observe that the choice of $P$ in our driving examples
(Subsections \ref{SSE:HEATEQUATION}-\ref{SSE:DELAYEQUATION})
will consider cases
the case when $\Imm P$ is closed and finite dimensional.
It is then useful to recall that, when the image of $P$ is closed, we can identify the space $B_b^P(H)$ with
$B_b(Im P)$ (and the same for the others). In particular,
in the case of Subsection \ref{SSE:DELAYEQUATION},
when $Im P=\R^n\times \{0\}$,
we immediately see that
$B^{P}_b(H) \sim B_b(\R^n)$, $C^{P}_b(H) \sim C_b(\R^n)$,
$UC^{P}_b(H) \sim UC_b(\R^n)$. This will be used in the sequel.
\hfill\qedo
\end{remark}

\noindent To prove our partial smoothing result we need
the following controllability-like assumption.
\begin{hypothesis}\label{ip:NC}
\begin{itemize}
\item[]
\item [(i)]
We have
\begin{equation}\label{eq:inclusionsmoothingC}
\operatorname{Im}\overline{Pe^{tA}}C\subseteq
\operatorname{Im}(P Q_t P^*)^{1/2},\qquad \forall t>0;
\end{equation}
Consequently, by the Closed Graph Theorem, the operator
$$
\Lambda^{P,C}(t):K\to H, \qquad
\Lambda^{P,C}(t)k:=(P Q_t P^*)^{-1/2}\overline{Pe^{tA}}Ck
\quad \forall k \in K,
$$
is well defined and bounded for all $t>0$.
\item [(ii)]
For every $T>0$ there exists $\kappa_T>0$ and $\gamma \in (0,1)$ such that
$$
\|\Lambda^{P,C}(t)\|_{\call(K,H)} \le \kappa_T t^{-\gamma}, \qquad
\forall t \in (0,T].
$$
\end{itemize}
\end{hypothesis}
Hypothesis \ref{ip:NC}-(i) is the analogous of the null controllability assumption which guarantees the strong Feller property of the associated
Ornstein-Uhlenbeck transition semigroup, see e.g. \cite{DP1} and \cite{Z}, while \ref{ip:NC}-(ii) is an asumption that guarantees that for $t\to 0$, the operator norm of $\Lambda^{P,C}(t)$ blows up in an integrable way. Both the assumptions can be verified in some models, namely in the following we show that the motivating examples introduced in Section \ref{SE:EXCONTROL} satisfy Hypothesis \ref{ip:NC}.
\begin{remark}\label{rm:NCexamples}
In the case of Subsection \ref{SSE:HEATEQUATION}
Hypothesis \ref{ip:NC} is satisfied, e.g., if we choose:
\begin{itemize}
  \item $H, \overline{H},A,C$ as in
Remark \ref{rm:computeexamples},
  \item $Q=(-A)^{-2\beta}$ for some $\beta\ge 0$
  \item $P$ a projection on a finite dimensional subspace contained in $\cald(-A)^{-\alpha}$ for some
      $\alpha>\beta+ \frac14$.
\end{itemize}
See Appendix  A1.
\smallskip
\\
In the case of Subsection \ref{SSE:DELAYEQUATION}
the above Hypothesis \ref{ip:NC}
are satisfied if:
\begin{itemize}
  \item we choose $H, \overline{H},A,C,P$ as in
Remark \ref{rm:computeexamples};
  \item we assume that {$
\operatorname{Im}\left(e^{ta_0}b_0 +\dis\int_{-d}^0 1_{[-t,0]}e^{(t+r)a_0}b_1(dr)
\right)
\subseteq\operatorname{Im}\sigma,
\quad \forall t>0.
$} We notice that this condition is verified when $\sigma$ is invertible, and it is a weaker assumption.
\end{itemize}
See Appendix  A.2.
\hfill\qedo
\end{remark}

Now we give the result.

\begin{proposition}\label{prop:partsmooth}
Let Hypotheses \ref{ip-sde-common}, \ref{ip:PC} and \ref{ip:NC}-(i) hold true.
\\
Then the semigroup $R_t,\,t>0$ maps functions $\phi\in B_b^P(H)$
into functions which are $C$-Fr\'echet differentiable in $\overline{H}$, and the $C$-derivative is given, for all $x \in \overline{H}$, by
\begin{align}\label{eq:formulader-gen-P}
 \nabla^C(R_{t}[\phi])(x)k &=\int_{H}\bar\phi\left(z_1+Pe^{tA}x\right)
\<\Lambda^{P,C}(t) k,
(PQ_tP^*)^{-1/2}z_1\>_H\caln(0,PQ_tP^*)(dz_1)
\\
&=
\E\left[\bar\phi\left(PX(t;x)\right)
\<\Lambda^{P,C}(t) k,
(PQ_tP^*)^{-1/2}PW_A(t)\>
\right]
\end{align}
Moreover, for any $\phi\in B^P_b(H)$ and any $k\in K$,
\begin{equation}\label{norm-Cder}
\vert \<\nabla^C R_t[\phi](x), k\>\vert \leq
\Vert \Lambda^{P,C}(t) \Vert_{\call(K, H)} \Vert \phi\Vert_\infty \vert k\vert.
\end{equation}
Furthermore, if $\phi\in C^P_b(H)$, then $\nabla^C R_t[\phi]\in C((0,T]\times \overline{H};K)$.
Finally, if also Hypothesis \ref{ip:NC}-(ii) holds, then the map $(t,x)\to R_t[\phi](x)$ belongs to $C_\gamma^{0,1,C}([0,T]\times \overline{H})$.
\end{proposition}

\dim
\\
If $\phi \in B_b^P(H)$, then, by \eqref{ornstein-sem-gen},
for every $t>0$ and $x \in H$,
\begin{equation}
 \label{eq:ornstein-sem-phibarCV}
R_t[\phi](x)
=\int_H \bar\phi(Pz+Pe^{tA}x)\caln(0,Q_t)(dz),
=\int_{H} \bar\phi(z_1+Pe^{tA}x)\caln(0,PQ_tP^*)(dz_1),
\end{equation}
where we adopt the change of variable $z_1=P z$ and we used that the image of the measure $\caln(0,Q_t)$ through
$P:H\to H$ is, clearly, $\caln(0,PQ_tP^*)$.
Now notice that, defining, for $t>0$,
$$
\call_t:L^2(0,t;K)\to H, \qquad \call_t u = \int_0^t e^{(t-s)A}Gu(s)ds,
$$
we get, by simple computations, that
$$
|(P\call_t)^*x|^2=\<Q_tP^*x,P^*x\>
\quad \hbox{which implies}\quad
\operatorname{Im} P\call_t= \operatorname{Im} (PQ_t P^*)^{1/2}.
$$
Hence, in particular the image of of $(PQ_t P^*)^{1/2}$ is contained in $\Imm P$. Moreover, if Hypothesis \ref{ip:NC}-(i) holds, the above also implies that
$\Imm\overline{Pe^{tA}}C\subseteq  \Imm P$ for all $t>0$.
Using this fact, for $\phi \in B^P_b(H)$,
$t>0$, $x \in H$, $k\in K$,
$\alpha \in \R$,
\begin{align}
R_t[\phi](x+\alpha Ck)=&\int_H
\bar\phi(Pz+\overline{Pe^{tA}}(x+\alpha Ck)\caln(0,Q_t)(dz),
\nonumber
\\
=&\int_{H} \bar\phi(z_1+\overline{Pe^{tA}}(x+\alpha Ck)
\caln(0,PQ_t P^*)(dz_1),
\label{eq:ornstein-sem-phibarCVnew}
\end{align}
where, in the second equality, we still use
the change of variable $z_1= Pz$.
Now we apply the change of variable
$z_2=z_1+\overline{Pe^{tA}}\alpha Ck$ to
\eqref{eq:ornstein-sem-phibarCVnew} getting that,
for every $t>0$ and $\phi \in B_b^P(H)$,
\begin{align}
&R_t[\phi](x+\alpha Ck)
=\int_{H} \bar\phi(z_2+Pe^{tA}x)
\caln(\alpha\overline{Pe^{tA}} Ck,PQ_tP^*)(dz_2).
 \label{eq:ornstein-sem-phibarCVnewbis}
\end{align}
Now, for $\phi\in B^P_b(H)$, $x\in H$, $k\in K$, $\alpha \in \R-\{0\}$, we get, by \eqref{eq:ornstein-sem-phibarCV}-\eqref{eq:ornstein-sem-phibarCVnewbis},
\begin{align}
\label{eq:rappincr}
\frac{1}{\alpha}
&\left[R_{t}[\phi](x+\alpha Ck)-R_{t}[\phi](x)\right]=
\\[2mm]
\notag
=&\frac{1}{\alpha}
\left[\int_{H} \bar\phi(z_1+Pe^{tA}x)
\caln(\alpha\overline{Pe^{tA}} Ck,PQ_tP^*)(dz_1)
-\int_{H}\bar\phi\left(z_1+Pe^{tA}x\right)
\caln\left(0,PQ_{t}P^*\right)(dz_1)\right].
\end{align}
By the Cameron-Martin theorem, see
e.g. \cite{DP3}, Theorem 1.3.6, the Gaussian measures
$\caln\left(\alpha\overline{Pe^{tA}} Ck,PQ_tP^*\right)$ and
$\mathcal{N}\left(0,PQ_tP^*\right)$ are equivalent if and only if
$\overline{Pe^{tA}}Ck\in\operatorname{Im}(PQ_tP^*)^{1/2}$.
In such case, setting, for $y \in \operatorname{Im}(PQ_tP^*)^{1/2}$, the density is
\begin{align}
&d(t,y,z):=\frac{d\caln\left(y,PQ_tP^*\right)}
{d\mathcal{N}\left(0,PQ_tP^*\right)  }(z)
\nonumber \\
&  =\exp\left\{  \left\langle (PQ_tP^*)^{-1/2}
y,(PQ_tP^*)^{-1/2}z\right\rangle_H
-\frac{1}{2}\left|(PQ_tP^*)^{-1/2}y\right|_H^{2}\right\} .
\label{eq:density1}
\end{align}
Such density is well defined for $z\in (\ker PQ_tP^*)^\perp$
(see e.g. \cite[Proposition 1.59]{FabbriGozziSwiech}).
Hence, by \eqref{eq:rappincr},
\begin{multline}
\label{eq:incrnew}\lim_{\alpha\rightarrow 0}\frac{1}{\alpha}
\left[R_{t}[\phi](x+\alpha Ck)-R_{t}[\phi](x)\right]=
\\
\lim_{\alpha\rightarrow 0}
\int_{H}\bar\phi\left(z_1+Pe^{tA}x\right)
\frac{d(t,\alpha \overline{Pe^{tA}}Ck,z_1)-1}{\alpha}
\caln(0,PQ_tP^*)(dz_1)
\end{multline}
Now we observe that, by the definition of $\Lambda^{P,C}(t)$,
$$
\frac{d(t,\alpha \overline{Pe^{tA}}Ck,z_1)-1}{\alpha}
=\frac{1}{\alpha}\left[\exp\left\{\alpha
\left\langle
\Lambda^{P,C}(t) k,(PQ_tP^*)^{-1/2}z_1\right\rangle_H
-\frac{\alpha^2}{2}
\left|\Lambda^{P,C}(t)k\right|_H^{2}\right\}
-1\right].
$$
When $\alpha \to 0$ the above limit is, $\left\langle\Lambda^{P,C}(t)k,(PQ_tP^*)^{-1/2}z_1\right\rangle_H$, which makes sense for all $z_1\in (\ker PQ_tP^*)^\perp$
and is an $L^2(H;\caln(0,PQ_tP^*))$ function of $z_1$
(see again, e.g., \cite[Proposition 1.59]{FabbriGozziSwiech}).
Moreover, with respect to the measure $\caln(0,PQ_tP^*)(dz_1)$ the map
$$
z_1 \mapsto \mathcal{Q}_t (z_1):= \left\langle\Lambda^{P,C}(t)k,(PQ_tP^*)^{-1/2}z_1\right\rangle_H
$$
is real valued Gaussian random variable with mean $0$ and
variance $\left|\Lambda^{P,C}(t)k\right|_H^{2}$ (see, e.g. \cite[Remark 2.2]{DPZ91}).
So in particular, for all $L>0$ $\E[e^{L|\mathcal{Q}_t| }]<+\infty$.
Now it is easy to see that
$$
\frac{d(t,\alpha \overline{Pe^{tA}}Ck,z_1)-1}{\alpha}
\le
e^{|\mathcal{Q}_t|+\left|\Lambda^{P,C}(t)k\right|_H^{2}}
$$
Hence we can apply the dominated convergence theorem
to \eqref{eq:incrnew} getting
\begin{align*}
& \exists \lim_{\alpha\rightarrow 0}\frac{1}{\alpha}
\left[R_{t}[\phi](x+\alpha Ck)-R_{t}[\phi](x)\right]=
\\
&\lim_{\alpha\rightarrow 0}\frac{1}{\alpha}
\int_{H}\bar\phi\left(z_1+Pe^{tA}x\right)
\frac{d(t,\alpha\overline{Pe^{tA}} Ck,z_1)-1}{\alpha}
\caln(0,PQ_tP^*)(dz_1)
\\[2mm]
&  =\int_{H}\bar\phi\left(z_1+Pe^{tA}x\right)
\<\Lambda^{P,C}(t) k, (PQ_tP^*)^{-1/2}z_1\>_H
\caln(0,PQ_tP^*)(dz_1)
\end{align*}
Consequently, along Definition \ref{df4:Gderunbounded}-(i),
there exists the $C$-directional derivative
$\nabla^C R_{t}\left[\phi\right](x;k)$
which is equal to the above right hand side.
Using that $\Lambda^{P,C}(t)$ is continuous we
see that the above limit is uniform for $k$ in the unit ball of $K$, so there exists the $C$-Fr\'echet derivative
$\nabla^C R_{t}\left[\phi\right](x)$.
From the above and from
\cite[Proposition 1.59]{FabbriGozziSwiech} we get
\begin{align*}
|\nabla^C R_{t}\left[\phi\right](x;k)|
&\leq
\Vert \bar\phi\Vert_\infty
\left(\int_{H}
\<\Lambda^{P,C}(t) k, (PQ_tP^*)^{-1/2}z_1\>_H^2
\caln(0,PQ_tP^*)(dz_1)\right)^{1/2}
\\[3mm]
\nonumber
&
= \Vert \bar\phi\Vert_\infty
\Vert \Lambda^{P,C}(t)k  \Vert_{H}
\le \Vert\phi\Vert_\infty
\Vert \Lambda^{P,C}(t)\Vert_{\call(K;H)} |k|_K . \nonumber
\end{align*}
This gives the required estimate.
The statement on continuity follows
using the same argument as in \cite[Theorem 4.41-(ii)]{FabbriGozziSwiech}.
The last statement follows by the last part of Definition \ref{df4:Gspaces}.
\qed


\begin{remark}\label{rm:smoothingproof}
The proof generalizes the one of Theorem 4.1 in \cite{FGFM-I} and the one of Theorem 4.41 in \cite{FabbriGozziSwiech}. The main difference between Theorem 4.1 in \cite{FGFM-I} and the present proposition is that here we are able to handle an unbounded operator $C$ by enlarging the space $H$. Notice that in the proof $C$ appears only through the operator $\overline{Pe^{tA}}C$. Notice also that, as proved above, the image of such operator is contained in $\operatorname{Im} P$, which is not obvious due to the presence of the closure.
On the other hand, the difference with respect Theorem 4.41 in \cite{FabbriGozziSwiech} is that there $P$ is missing and the
partial derivatives are taken in the unbounded but less general case of Definition \cite[Definition 4.4]{FabbriGozziSwiech}.
\hfill\qedo
\end{remark}

\begin{remark}\label{rm:partsmooth-second}
Generalizing to our setting the ideas of Proposition 4.5 in \cite{FGFM-I}
it is possible to prove that, if $\phi$ is more regular
(i.e. $\phi\in C^1_b(H)\cap C_b^P(H)$,
also $\nabla^C R_t[\phi]$ has more regularity, i.e.
$\nabla\nabla^{C}R_{t}\left[\phi\right]$,
$\nabla^{C}\nabla R_{t}\left[\phi\right]$ exist, coincide, and
satisfy suitable formulae and estimates.
We omit them here since we do not need them for the purpose of this paper. They will be useful to find optimal feedback controls, which will be the subject of a subsequent paper.
\hfill\qedo\end{remark}

\section{Partial smoothing for convolutions}
\label{sec:convpartsmooth-abstr-setting}

To solve HJB equations like \eqref{HJBformale-diri} and \eqref{HJBformale-delay}
we need to extend the partial smoothing result of the previous section to convolutions.

We need first to introduce suitable spaces where such convolutions live
and which will be useful later to perform the fixed point argument to find the solution of our HJB equations.

\begin{definition}\label{df:Sigma}
Let $T>0$, $\eta \in (0,1)$.
A function $g\in C_b([0,T]\times H)\cap
C_b((0,T]\times \overline{H})$
belongs to $\Sigma^1_{T,\eta}$ if
\begin{itemize}
  \item
there exists a function
$f\in C_b([0,T]\times H)$ such that\footnote{By continuity this also implies $g(0,x)=f(0,Px)$ for all $x \in H$.}
$$g(t,x)=f\left(t,\overline{Pe^{tA}}x\right),
\qquad \forall (t,x) \in (0,T]\times \overline{H};
$$
  \item
for any $t\in(0,T]$ the function $g(t,\cdot)$ is
$C$-Fr\'echet differentiable on $\overline{H}$ and there exists
a function $\bar f\in C_b((0,T]\times H;K)$
such that
$$
t^\eta \nabla^C g(t,x)=\bar f\left(t,\overline{Pe^{tA}}x\right),
\qquad \forall (t,x) \in (0,T]\times \overline{H}.
$$
\end{itemize}
\end{definition}

\begin{remark}\label{rm:partsmooth-second-1}
Arguing as in \cite[Section 5]{FGFM-I}, it is possible to define a subspace of $\Sigma^1_{T,\eta}$ of functions $g$ such that there exists the second order derivative $\nabla\nabla^C$ which depends in a special way on $x\in \overline{H}$. This could be useful to prove second order regularity of  the solution of our HJB equations. As in Remark \ref{rm:partsmooth-second} we omit this step here: it will be useful to find optimal feedback controls, which will be the subject of a subsequent paper.
\hfill\qedo\end{remark}

\begin{remark}
\label{rm:sigma1}
We observe that, using substantially the same argument as
\cite[Lemma 5.2]{FGFM-I}, one can prove that
\begin{equation}\label{eq:SigmaSubspace}
  \Sigma^1_{T,\eta} \hbox{ is a closed subspace of }
  C_\eta^{0,1,C}([0,T]\times \overline{H})
\end{equation}
Moreover we also observe that, by Proposition \ref{prop:partsmooth},
it is immediate to see that, under our Hypotheses \ref{ip-sde-common}-\ref{ip:PC}-\ref{ip:NC}, for any $\phi \in C_b^P(H)$, we have
$R_t[\phi] \in \Sigma^1_{T,\gamma}$.
\end{remark}

We now come back to the abstract common setting and we state a first lemma on the regularity of the convolution type terms.
\begin{lemma}\label{lemma_convoluzione}
Let Hypotheses \ref{ip-sde-common}, \ref{ip:PC} and \ref{ip:NC} hold.
Let $T>0$, $C\in\call(K,\overline{H})$ and let $\psi:K^*\rightarrow\R$ be a Lipschitz continuous function.
For every $g \in \Sigma^1_{T,\gamma}$ (where $\gamma$ is
given in Hypothesis \ref{ip:NC}-(ii)), the function
$\hat{g}:[0,T]\times \overline{H} \rightarrow \R$ belongs to
$\Sigma^1_{T,\gamma}$ where
\begin{equation}\label{iterata-primag}
\hat{g}(t,x) =\int_{0}^{t}
R_{t-s} [\psi(\nabla^{C}g(s,\cdot))](x)  ds, \qquad (t,x) \in
[0,T]\times \overline{H}.
\end{equation}
Hence, in particular, $\hat g(t,\cdot)$
is $C$-Fr\'echet differentiable on $\overline{H}$ for every $t\in (0,T]$ and, for a suitable constant $\hat\kappa$ (depending only on $T$ and $\psi$),
\begin{equation}\label{stimaiterata-primag}
\left\vert \nabla ^C(\hat g(t,\cdot))(x) \right\vert_{K^*}    \leq
    \hat\kappa\left(t^{{1-\gamma}}+t^{{1-2\gamma}}
    \Vert g \Vert_{C^{0,1,C}_{{\gamma}}}\right),
\qquad \forall (t,x)\in (0,T]\times \overline{H}.
 \end{equation}
Moreover, for every $g_1,g_2 \in \Sigma^1_{T,\gamma}$ (where $\gamma$ is given in Hypothesis \ref{ip:NC}-(ii)), the function
$\hat{g_1}-\hat{g_2}:[0,T]\times \overline{H} \rightarrow \R$ belongs to $\Sigma^1_{T,\gamma}$
and, for a suitable constant $\kappa$ (depending only on $T$ and $\psi$),
\begin{multline}\label{stimaiterata-g1g2}
\left\vert
\hat g_1(t,x)
-\hat g_2(t,x))
 \right\vert +  t^\gamma \left\vert
\nabla ^C(\hat g_1(t,\cdot))(x)
-\nabla ^C(\hat g_2(t,\cdot))(x)
 \right\vert_{K^*}
 \\
 \leq
    \kappa\left(t+t^{{1-\gamma}}\right)
    \Vert g_1-g_2 \Vert_{C^{0,1,C}_{{\gamma}}},
\qquad \forall (t,x)\in (0,T]\times \overline{H}.
 \end{multline}
%
%
\end{lemma}
\dim
We start by proving that $\hat g$ from (\ref{iterata-primag}) is well defined, continuous, and $C$-Fr\'echet differentiable.
First of all, for any $g\in\Sigma^1_{T,\gamma}$, we denote by $f_g$ and $\bar f_g$ the functions associated to it in Definition \ref{df:Sigma}.
Hence, given any $g\in\Sigma^1_{T,\gamma}$, we have, for
$0<s\le t$, $x,z\in H$:
\begin{align}\label{eq:nablaCshiftg}
 & \psi\left(\nabla^{C}g(s,z+e^{(t-s)A}x)\right)
 =
\psi\left( s^{-\gamma}\bar f_g\left(s, \overline{Pe^{sA}}z+\overline{Pe^{tA}}x\right)\right)
\end{align}
Hence we can give meaning to the left hand side also for
$x \in \overline{H}$. So we can write
\begin{align}
\notag
&  \int_{0}^{t}
R_{t-s} \left[\psi\left(\nabla^{C}(g(s,\cdot))\right)\right](x)  ds
=\int_{0}^{t}
\int_H \psi\left(\nabla^{C}g(s,z+e^{(t-s)A}x)\right)
\caln(0,Q_{t-s})(dz)
\\
&=\int_{0}^{t}
\int_H \psi\left( s^{-\gamma}\bar f_g\left(s, {Pe^{sA}}z+\overline{Pe^{tA}}x\right)\right)
\caln(0,Q_{t-s})(dz),
\qquad \forall(t,x)\in [0,T]\times\overline{H}.
\label{eq:defRconv1}
\end{align}
The above implies that $\hat g$ is well defined on $[0,T]\times\overline{H}$. Continuity follows
using the same argument as in
\cite[Proposition 4.50-(ii)]{FabbriGozziSwiech}.
Consequently the function $f_{\hat g}$ associated to $\hat g$ along Definition \ref{df:Sigma} is
\begin{align*}
f_{\hat g} (t,y)&=\int_{0}^{t}
\int_{H} \psi\left(s^{-\gamma}
 \bar f_g\left(s, {Pe^{sA}}z+y\right)
\right)\caln(0,Q_{t-s})(dz)
\end{align*}
and, by the Lipschitz assumptions on $\psi$,
$$
\Vert\hat f_{\hat g} \Vert_\infty \le \int_{0}^{T}[\psi]_{Lip}
\left(|\psi(0)| + s^{-\gamma} \Vert \bar f_g \Vert_\infty \right)ds
\le
[\psi]_{Lip}\left[|\psi(0)| T+ \Vert \bar f_g \Vert_\infty (1-\gamma)^{-1} T^{1-\gamma}\right]
$$
To compute the $C$-derivative we first compute, using what is given above,
\begin{align}
\label{eq:defRconv2}
&\int_{0}^{t}
 R_{t-s}\left[\psi\left( \nabla^{C}(g(s,\cdot))\right)\right](x+\alpha Ck)ds \\
 \notag
&=\int_{0}^{t} \int_H  \psi\left( s^{-\gamma}
\bar f_g\left(s, {Pe^{sA}}z+\overline{Pe^{tA}}(x+\alpha Ch)\right)
\right)\caln(0,Q_{t-s})(dz) ds
\\[2mm]
\notag&=\int_{0}^{t} \int_H
\psi\left(s^{-\gamma}
\bar f_g\left(s, Pe^{sA}z+\overline{Pe^{tA}}x\right)
\right)
\caln\left(
\overline{Pe^{tA}}\alpha Ck,Q_{t-s}\right)(dz)ds
\\[2mm]
\notag
&=\int_{0}^{t} \int_H
\psi\left(s^{-\gamma}
\bar f_g\left(s, Pe^{sA}z+\overline{Pe^{tA}}x\right)
\right)
d(t-s,\alpha \overline{Pe^{tA}}Ck,z)
\caln\left(0,Q_{t-s}\right)(dz)ds,
\end{align}
where, in the last two equalities, we have used Cameron-Martin Theorem
as in the proof of the above Proposition \ref{prop:partsmooth}, and
the density $d$ is given by \eqref{eq:density1}.
Hence, using \eqref{eq:defRconv1}-\eqref{eq:defRconv2},
\begin{align*}
&\lim_{\alpha\rightarrow 0}\dfrac{1}{\alpha}
 \left[\int_{0}^{t}
 R_{t-s}\left[\psi\left( \nabla^{C}(g(s,\cdot))\right)\right](x+\alpha Ck)ds -
\int_{0}^{t}
 R_{t-s}\left[\psi\left( \nabla^{C}(g(s,\cdot))\right)\right](x)  ds\right]=
 \\[2mm]
 &=\lim_{\alpha\rightarrow 0}\dfrac{1}{\alpha}
 \int_{0}^{t} \int_H
 \psi\left(s^{-{\gamma}}
 \bar f_g\left(s, Pe^{sA}z+ \overline{Pe^{tA}}x\right)
 \right)
 \frac{d(t-s,\alpha  \overline{Pe^{tA}}Ck,z)-1}{\alpha}
 \caln\left(0,Q_{t-s}\right)(dz)ds
\end{align*}
At this point we argue exactly as in the proof of the above Proposition \ref{prop:partsmooth} getting, uniformly for $k$ in the unit sphere of $K$:
\begin{align*}
&\lim_{\alpha\rightarrow 0}\dfrac{1}{\alpha}
 \left[\int_{0}^{t}
 R_{t-s}\left[\psi\left( \nabla^{C}(g(s,\cdot))\right)\right](x+\alpha Ck)ds -
\int_{0}^{t}
 R_{t-s}\left[\psi\left( \nabla^{C}(g(s,\cdot))\right)\right](x)  ds\right]=
 \\[2mm]
&=
 \int_{0}^{t} \int_H
 \psi\left(s^{-\gamma}
 \bar f_g\left(s, Pe^{sA}z+\overline{Pe^{tA}}x\right)\right)
 \<Q_{t-s}^{-1/2}
 \overline{Pe^{tA}} Ck, Q_{t-s}^{-1/2}z\>_H
 \caln\left(0,Q_{t-s}\right)(dz)ds.
\end{align*}
This implies the required $C$-Fr\'echet differentiability and
\begin{align}
&\<\nabla ^C \hat g(t,x),k \>_{K}=
\label{eq:derCConv}
\\
\nonumber
&=
\int_{0}^{t} \int_H
\psi\left(
s^{-\gamma} \bar f_g\left(s, Pe^{sA}z+\overline{Pe^{tA}}x\right)
\right)
\<(Q_{t-s}^{-{1/2}}
\overline{Pe^{tA}}Ck, Q_{t-s}^{-{1/2}}Pz\>_{H}
\caln\left(0,Q_{t-s}\right)(dz)ds.
\end{align}
Moreover, the right hand side of \eqref{eq:derCConv} provides, when we substitute $\overline{Pe^{tA}}x$ with $y$, the function $\bar{{f}}_{\hat g}$
associated to $\hat g$ along the second part of Definition \ref{df:Sigma}.
\\
At this point, in order to prove estimate (\ref{stimaiterata-primag}),
we use the above representation and the Holder inequality:
\begin{align*}
&\left\vert\<\nabla ^C \hat g(t,x),k \>_K\right\vert \le
\\
 &\leq [\psi]_{Lip}\int_{0}^{t}
  \int_H  \left(|\psi(0)|+\left\vert
 s^{-{\gamma}} \bar f_g\left(s, Pe^{sA}z+\overline{Pe^{tA}}x\right)
 \right\vert\right)\\
 &\quad\left \vert
 \<(Q_{t-s})^{-{1/2}}\overline{Pe^{tA}} Ck, Q_{t-s})^{-{1/2}}Pz\>_{H}
  \right\vert\caln(0,Q_{t-s})(dz) ds
 \\
 &\leq [\psi]_{Lip}\int_{0}^{t}\left(|\psi(0)|+s^{-\gamma} \left\Vert g \right\Vert_{C^{0,1,C}_{\gamma}}
 \right)
 \left\Vert Q_{t-s}^{-{1/2}}\overline{Pe^{tA}}Ck \right\Vert_{\call(K;H)} ds \\
  &\leq  \kappa_T [\psi]_{Lip}\int_{0}^{t}\left(|\psi(0)|+s^{-{\gamma}}
  \left\Vert g \right\Vert_{C^{0,1,C}_{{\gamma}}}\right)
  (t-s)^{-\gamma}\vert k\vert_{K} \,ds
\end{align*}
Since
$$
\int_{0}^{t}|\psi(0)|
  (t-s)^{-\gamma}\vert k\vert_{K} \,ds
= |\psi(0)|\vert k\vert_{K}
\frac{1}{1-\gamma}t^{1-\gamma}
$$
$$
\int_{0}^{t}s^{-{\gamma}}
  \left\Vert g \right\Vert_{C^{0,1,C}_{{\gamma}}}
  (t-s)^{-\gamma}\vert k\vert_{K} \,ds
= \left\Vert g \right\Vert_{C^{0,1,C}_{{\gamma}}}\vert k\vert_{K}
\int_{0}^{t}s^{-{\gamma}}(t-s)^{-\gamma} \,ds
=
\left\Vert g \right\Vert_{C^{0,1,C}_{{\gamma}}}\vert k\vert_{K}
\beta(1-\gamma,1-\gamma)t^{1-2\gamma},
$$
where by $\beta(\cdot,\cdot)$ we mean the Euler Beta function, the claim follows.
\\
The proof of \eqref{stimaiterata-g1g2} follows in a similar way, taking into account estimate \eqref{stimaiterata-primag} an the fact that $\psi$ is Lipschitz continuous:
\begin{multline*}\label{stimaiterata-g1g2}
\left\vert
\hat g_1(t,x)
-\hat g_2(t,x))
 \right\vert +  t^\gamma \left\vert
\nabla ^C(\hat g_1(t,\cdot))(x)
-\nabla ^C(\hat g_2(t,\cdot))(x)
 \right\vert_{K^*}
 \\
\leq \vert  \int_{0}^{t}
R_{t-s} \left[\psi\left(\nabla^{C}(g_1(s,\cdot))\right)-
\psi\left(\nabla^{C}(g_2(s,\cdot))\right)\right](x) ds\vert\\
+t^\gamma\vert \nabla^{C} \int_{0}^{t}
R_{t-s} \left[\psi\left(\nabla^{C}(g_1(s,\cdot))\right)
-\psi\left(\nabla^{C}(g_2(s,\cdot))\right)\right](x) ds\vert\\
\leq
t^\gamma\left\vert\int_{0}^{t} \int_H \left(
\psi\left(
s^{-\gamma} \bar f_{g_1}\left(s, Pe^{sA}z+\overline{Pe^{tA}}x\right)
\right)-\psi\left(
s^{-\gamma} \bar f_{g_2}\left(s, Pe^{sA}z+\overline{Pe^{tA}}x\right)
\right)\right)
\right.\\
\left.\<(Q_{t-s}^{-{1/2}}
\overline{Pe^{tA}}Ck, Q_{t-s}^{-{1/2}}Pz\>_{H}
\caln\left(0,Q_{t-s}\right)(dz)ds\right\vert + \kappa t \Vert g_1-g_2 \Vert_{C^{0,1,C}_{{\gamma}}}\\
\leq
    \kappa\left(t+t^{{1-\gamma}}\right)
    \Vert g_1-g_2 \Vert_{C^{0,1,C}_{{\gamma}}},
\qquad \forall (t,x)\in (0,T]\times \overline{H}.
 \end{multline*}

 \qed

\section{Applying partial smoothing to stochastic control problems}
\label{sec-HJB}

In this Section we first present the stochastic optimal control problem we aim to treat, then we show how the theory developed in the previous sections allows us to solve the associated HJB equation.

\subsection{A stochastic control problem in the abstract setting}
\label{subsec-contr.pr-abstract}
We consider the setting of Hypotheses \ref{ip-sde-common}-\ref{ip:PC}-\ref{ip:NC} which we assume to hold.
We present first the objective functional and then the state equation.
The goal of the controller is to minimize the following
finite horizon cost (here $X(\cdot;t,x)$ is the state process starting at time $t$ with the datum $x$),
\begin{equation}\label{costoastratto-common}
J(t,x;u)=\E \left(\int_t^T \left[\ell_0(s)+\ell_1(u(s))\right]\,ds + \phi(X(T;t,x))\right)
\end{equation}
over all controls $u(\cdot)$ in
\begin{equation}\label{eq:admcontr}
  \calu:=\left\{
  u:[0,T]\times \Omega \to U \subseteq K, \text{ progressively measurable}
  \right\}
\end{equation}
under the following assumption.
\begin{hypothesis}\label{ip-costo}
We assume that:
\begin{itemize}
  \item[(i)]
the final cost $\phi$ belongs to $C_b^P(H)$ (see Definition \ref{df:spaziphi1});
\item[(ii)]
The current cost $\ell_0$ is measurable and bounded;
  \item[(iii)]
the set $U\subset K$ is closed and bounded and the current cost $\ell_1:U\rightarrow \R$ is measurable and bounded from below.
\end{itemize}
\end{hypothesis}

\begin{remark}
Note that here the current cost does not depend on the state. Putting the dependence on the state in the current cost would increase considerably the technical arguments, namely the fixed point argument in the proof of Theorem \ref{esistenzaHJB}, Section \ref{sub:HJBsol}, wouldn't work, and it is left for a further paper.
\newline We also underline that the above technical problem cannot be overcome by transforming our Bolza type problem in a Mayer type problem (i.e. a problem with only the terminal cost) on the line of what is done e.g. in 
\cite[Remark 7.4.1, p.714]{CannarsaSinestrari04}.
Indeed, using such transformation the state dependent running cost would disappear but the state equation would become nonlinear. This would prevent the use of our results on partial smoothing which, up to now, apply only to linear state equations.
\end{remark}

Before introducing the state equation we observe that,
due to Hypothesis \ref{ip-costo}-(i) above, what matters for the controller
is the process $PX(\cdot)$.
Now consider the following controlled SDE in the real separable Hilbert space $H$ (here $0\le t \le s \le T$)
\begin{equation}\label{eq-common-contr}
  \left\{
  \begin{array}{l}
  \dis
d X(s)= [AX(s)+Cu(s)]ds +Q^{1/2}dW(s), \qquad s\in (t,T],
\\\dis
X(s)=x\in H, \quad \forall s \in [0,t]
\end{array}
\right.
\end{equation}
where $A,\,G$ and $W$ are as in Hypothesis \ref{ip-sde-common}
and $C\in\call(K,\overline{H})$ is as in Hypothesis \ref{ip:PC}.
Equation (\ref{eq-common-contr}) is formal and has to be considered in its mild formulation (using the so-called variation of constants, see e.g. \cite[Chapter 7]{DP1}) which still present some issues. Indeed the mild solution of \eqref{eq-common-contr} is, still formally,
\begin{equation}
X(s)=e^{(s-t)A}x+\int_t^s{e^{(s-r)A}}C u(r) dr +\int_t^se^{(s-r)A}Q^{1/2}dW(r),
\text{ \ \ \ }s\in[t,T]. \\
  \label{eq-mild-common}
\end{equation}
Here the first and the third term belong to $H$, thanks to Hypothesis \ref{ip-sde-common}, while the second in general does not.
By Hypothesis \ref{ip:PC}-(i) we see that the second term can be written as
$$
\int_t^s\overline{e^{(s-r)A}}Cu(r)dr \in \overline{H}.
$$
Hence, even when $x\in H$ the mild solution belongs to $\overline{H}$ but not to $H$. Moreover, still using Hypothesis \ref{ip:PC}-(i), we see that the mild solution makes sense for all $x\in \overline{H}$ and belongs to $\overline{H}$.
\\
On the other hand, thanks to Hypothesis \ref{ip:PC}-(iii), even when
$x\in \overline{H}$ the process $PX(s)$ belongs to $H$ and can be written as
\begin{equation}
PX(s)  =\overline{Pe^{(s-t)A}}x+\int_t^s\overline{Pe^{(s-r)A}}C u(r) dr +\int_t^s P e^{(s-r)A}Q^{1/2}dW(r),
\text{ \ \ \ }s\in(t,T]. \\
  \label{eq-mild-commonP}
\end{equation}
We define the value function related to this control problem, as usual, as
\begin{equation}\label{valuefunction-gen}
V(t,x):= \inf_{u \in \calu}J(t,x;u).
\end{equation}

\subsection{Solution of the HJB equation}
\label{sub:HJBsol}

We define the Hamiltonian as follows:
the current value Hamiltonian $H_{CV}$ is
$$
H_{CV}(p\,;u):=\<p,u\>_{K}+\ell_1(u)
$$
and the minimum value Hamiltonian is
\begin{equation}\label{psi1-gen}
H_{min}(p)=\inf_{u\in U}H_{CV}(p\,;u),
 \end{equation}
The HJB equation associated to the stochastic optimal control problem presented in the above Section \ref{subsec-contr.pr-abstract} is then, formally,
\begin{equation}\label{HJBformale-common}
  \left\{\begin{array}{l}\dis
-\frac{\partial v(t,x)}{\partial t}=\call [v(t,\cdot)](x) +\ell_0(t)+
H_{min} (\nabla^C v(t,x)),\qquad t\in [0,T],\,
x\in H,\\
\\
\dis v(T,x)=\phi(x),
\end{array}\right.
\end{equation}
Here the differential operator $\call$
is the infinitesimal generator of the transition
semigroup $(R_{t})_{0\leq t\leq T}$  defined in (\ref{ornstein-sem-gen}) related to the process $Z$ solution of equation
(\ref{ornstein-gen}), namely $\call$ is formally defined by
\begin{equation}\label{eq:ell-gen}
 \call[f](x)=\frac{1}{2} Tr \; Q\nabla^2 f(x)
+ \< x,A^*\nabla f(x)\>.
\end{equation}
\begin{definition}\label{defsolmildHJB}
We say that a
function $v:[0,T]\times H\rightarrow\R$ is a mild
solution of the HJB equation (\ref{HJBformale-common}) if the following
are satisfied for some $\gamma \in(0,1)$:
\begin{enumerate}

\item $v(T-\cdot, \cdot)\in C^{0,1,C}_{{\gamma}}\left([0,T]\times \overline{H}\right)$;

\item the integral equation
\begin{equation}\
v(t,x) =R_{T-t}[\phi](x)+\int_t^T R_{s-t}\left[
H_{min}(\nabla^C v(s,\cdot))+\ell_0(s)\right](x)\; ds,
\label{solmildHJB-common}
\end{equation}
is satisfied on $[0,T]\times H$.
\end{enumerate}
\end{definition}
We notice that the request in the previous definition \ref{defsolmildHJB}, point 1, implies that $\nabla^C_xv(t,x)$ can blow up like $(T-t)^{-\gamma}$.

We now prove existence and uniqueness of a mild solution of the HJB equation (\ref{HJBformale-diri}) and (\ref{HJBformale-delay}) by a fixed point argument.
\begin{remark}\label{rm:crescitapoli-HJB}
Since functions in  $C^{0,1,C}_{\gamma}\left([0,T]\times \overline{H}\right)$
are bounded (see Definition \ref{df4:Gspaces}),
the above Definition \ref{defsolmildHJB} requires, among other properties, that a mild solution is continuous and bounded up to $T$.
This constrains the assumptions on the data, e.g. it implies that the final datum $\phi$ must be continuous and bounded.
We may change this requirement in the above definition asking only measurability or only polynomial growth in $x$ so allowing for more general datum $\phi$ in Hypothesis \ref{ip-costo}-(i). Our main results will remain true with straightforward modifications (see \cite[Chapter 4]{FabbriGozziSwiech} for a treatment of such cases in the case of bounded control operators).
\\
Similarly we may weaken the request of Hypothesis \ref{ip-costo}-(iii) on the boundedness of the set $U$. This may result in the fact that the Hamiltonian $H_{min}$ is not Lipschitz continuous. This case, even if more difficult, could still be treated using the ideas of \cite{G2,FMloclip}.
\end{remark}

\begin{theorem}\label{esistenzaHJB}
Let Hypotheses \ref{ip-sde-common}, \ref{ip:PC}, \ref{ip:NC} and \ref{ip-costo} hold true.
Then the HJB equation (\ref{HJBformale-common})
admits a mild solution $v$ according to Definition \ref{defsolmildHJB}.
Moreover $v$ is unique among the functions $w$ such that $w(T-\cdot,\cdot)\in\Sigma_{T,\gamma}$
and it satisfies, for a suitable constant $\kappa_{1,T}>0$, the estimate
\begin{equation}\label{eq:stimavmainteo}
\Vert v(T-\cdot,\cdot)\Vert_{C^{0,1,C}_{{\gamma}}}\le \kappa_{1,T}\left(\Vert\phi \Vert_\infty
+\Vert\ell_0 \Vert_\infty \right).
\end{equation}

%
\end{theorem}
\dim
We first prove existence and uniqueness of a solution in $\Sigma^1_{T,{\gamma}}$, by using a fixed point argument in it. To this aim,
first we rewrite (\ref{solmildHJB-common}) in a forward way. Namely
if $v$ satisfies \myref{solmildHJB-common} then, setting $w(t,x):=v(T-t,x)$ for any
$(t,x)\in[0,T]\times H$, we get that $w$ satisfies
\begin{equation}
w(t,x) =R_{t}[\phi](x)+\int_0^t R_{t-s}[H_{min}(
\nabla^C w(s,\cdot))+\ell_0(s)](x)\; ds,\qquad t\in [0,T],\
x\in H,\label{solmildHJB-forward}
\end{equation}
which is the mild form of the forward HJB equation
\begin{equation}\label{HJBformaleforward-common}
  \left\{\begin{array}{l}\dis
\frac{\partial w(t,x)}{\partial t}=\call [w(t,\cdot)](x) +\ell_0(t)+
H_{min} (\nabla^C w(t,x)),\qquad t\in [0,T],\,
x\in H,\\
\\
\dis w(0,x)=\phi(x).
\end{array}\right.
\end{equation}
Referring to equation \eqref{solmildHJB-forward}, which is the mild version of \eqref{solmildHJB-common}, define the map $\Upsilon$ on $\Sigma^1_{T,\gamma}$ by setting, for $g\in \Sigma^1_{T,{\gamma}}$,
$$ \Upsilon[g](0,x):=\phi(x),$$
and, for $(t,x)\in (0,T]\times \overline{H}$,
\begin{equation}\label{Gamma}
 \Upsilon[g](t,x):=R_{t}[\phi](x)+
\int_0^t \ell_0(s) ds+
\int_0^t R_{t-s}[
H_{min}(\nabla^C g(s,\cdot))](x)\; ds.
\end{equation}
Using Proposition \ref{prop:partsmooth}, in particular
\eqref{eq:ornstein-sem-phibarCV}, \eqref{eq:formulader-gen-P}
and the last statement on continuity,
we see that the sum of the first two terms of \eqref{Gamma} belongs to $\Sigma^1_{T,\gamma}$ with
\begin{align*}
f(t,x)&=\int_H\bar\phi(z_1+x)\caln(0,PQ_tP^*)(dz_1)
+\int_0^t \ell_0(s) ds,\\
\bar f(t,x)&=\int_H\bar\phi(z_1+x)\<\Lambda^{P,C}(t),(PQ_tP^*)^{-1/2}z_1\>
\caln(0,PQ_tP^*)(dz_1),\qquad
\end{align*}
Moreover, we use Lemma \ref{lemma_convoluzione} (simply substituting the generic function $\psi$ with $H_{min}$) to deduce that the third term of \eqref{Gamma} belongs to $\Sigma^1_{T,\gamma}$.
Hence $\Upsilon$
is well defined in $\Sigma^1_{T,\gamma}$ with values in $\Sigma^1_{T,\gamma}$ itself.
\newline As stated in Remark \ref{rm:sigma1}, $\Sigma^1_{T,\gamma}$ is a closed subspace of $C^{0,1,C}_{{\gamma}}([0,T]\times \overline{H})$,
and so if $\Upsilon$ is a contraction,
by the Contraction Mapping Principle there exists a unique (in $\Sigma^1_{T,{\gamma}}$)
mild solution of (\ref{HJBformale-common}).
We then prove the contraction property of $\Upsilon$.

Let $g_1,g_2 \in \Sigma^1_{T,{\gamma}}$. We evaluate
$$\Vert \Upsilon (g_1)-\Upsilon (g_2)\Vert_{\Sigma^1_{T,\gamma}}=\Vert \Upsilon(g_1)-\Upsilon (g_2)\Vert_{C^{0,1,C}_{\gamma}}.$$
For every $(t,x)\in (0,T]\times \overline{H}$, we have
\begin{align*}
\Upsilon (g_1)(t,x)- \Upsilon(g_2)(t,x) =
\int_0^t R_{t-s}\left[H_{min}\left(\nabla^C g_1(s,\cdot)\right)
 -H_{min}\left(\nabla^C g_2(s,\cdot)\right)\right](x)ds
\end{align*}
Hence we can use the second part of Lemma
 \ref{lemma_convoluzione}, namely estimate \eqref{stimaiterata-g1g2}, to get
\begin{multline}\label{eq:stimaUpsilonNew}
|\Upsilon (g_1)(t,x)- \Upsilon(g_2)(t,x)|
+
t^{{\gamma}}\vert \nabla^C\Upsilon (g_1)(t,x) - \nabla^C\Upsilon(g_2)(t,x) \vert_{K^*}
\le
\kappa (t +t^{{1-\gamma}})\Vert g_1-g_2 \Vert_{C^{0,1,C}_{{\gamma}}}.
\end{multline}
Hence, if $T$ is sufficiently small, we get
that the map $\Upsilon$ is a contraction in $\Sigma^1_{T,\gamma}$
and, if we denote by $w$ its unique fixed point, then $v:=w(T-\cdot,\cdot)$
turns out to be a mild solution of the HJB equation (\ref{HJBformale-common}),
according to Definition \ref{defsolmildHJB}.

Since the constant $C$ is independent of $t$, the case of generic $T>0$ follows
by dividing the interval $[0,T]$
into a finite number of subintervals of length $\delta$ sufficiently small, or equivalently, as done in \cite{Mas},
by taking an equivalent norm with an adequate exponential weight, such as
\[
 \left\Vert f\right\Vert _{\eta,C^{0,1,C}_{\gamma}}
=\sup_{(t,x)\in(0,T]\times \overline{H}}
\vert e^{\eta t}f(t,x)\vert+
\sup_{(t,x)\in (0,T]\times \overline{H}}  e^{\eta t}t^{\gamma}
\left\Vert \nabla^C f\left(  t,x\right)
\right\Vert _{K^*}.
\]

\qed



\appendix

\begin{appendix}

\section{Appendix: verifying Hypothesis \ref{ip:NC} in our examples}

\subsection{The case of boundary control}\label{Sec-min-ene}
We consider the stochastic heat equation with boundary control \eqref{eqDiri}, reformulated as an abstract evolution equation \eqref{eqDiri-abstr-contr}, where we choose, according to Subsection \ref{SSE:HEATEQUATION} and Remark \ref{rm:computeexamples},
$$
H=L^2(\calo), \qquad
\overline H= \cald\left((-A_0)^{-3/4-\eps}\right)
=H^{-3/2-2\eps}(\calo)
\quad \hbox{(for suitable small $\eps>0$),}
$$
$A=A_0$, $C=B=(-A_0)D$ as from \eqref{notazioneB}.
Moreover, as from Remark \ref{rm:NCexamples} we take $Q=(-A)^{-2\beta}$ for some $\beta\ge 0$ and
$P$ a projection on a finite dimensional subspace contained in $(-A)^{-\alpha}$ for some $\alpha>\beta+ \frac14$.

The covariance operator $Q_t$ is given by
\begin{equation}\label{cov-heat}
Q_t=\int_0^t (-A_0)^{-2\beta}e^{2sA_0}\,ds=(-A_0)^{-2\beta-1}(I-e^{2tA_0}).
\end{equation}
Notice that it can be deduced by the strong Feller property of the heat transition semigroup that $\operatorname{Im}e^{tA}\subset \operatorname{Im}Q_t^{1/2}$, see e.g. \cite[Section 9.4 and Appendix B]{DP1} for a comprehensive bibliography. Now we estimate
$\Vert Q_t^{-1/2}e^{tA_0} (-A_0 D)\Vert $.

In the sequel we denote by $\lambda_k\geq0,\, k\geq 1,\; \lambda_k\nearrow +\infty,$ the opposite of the eigenvalues of the Laplace operator in $\calo$:
\[
A_0e_k=-\lambda_ke_k,\, k\geq 1.
\]
\begin{lemma}\label{Lemma Qt}
 Let $Q_t$ be defined in (\ref{cov-heat}). For every $\varepsilon\in (0,\dfrac{1}{4})$, we get, for some $C_0>0$,
 \begin{equation}\label{Q_t-norm}
  \Vert Q_t^{-1/2}e^{tA_0}(-A_0 D)\Vert
  \le C_0 t^{-\frac{5}{4}-\beta -\varepsilon} ,
 \end{equation}
\end{lemma}
\dim We notice that $ Q_t^{-1/2}e^{tA_0}(-A_0D)=Q_t^{-1/2}e^{tA_0}(-A_0)^{\frac{3}{4}+\varepsilon}D_\varepsilon)$, where $D_\varepsilon=(-A_0)^{\frac{1}{4}-\varepsilon}D$ is bounded $\forall\, \varepsilon \in \Big(0,\dfrac{1}{4}\Big)$.
Moreover for every $a\in \partial\,\calo$
\begin{align*}
\vert Q_t^{-1/2}e^{tA_0}(-A_0)^{\frac{3}{4}+\varepsilon}D_\varepsilon a\vert^2
&=\sum_{k=1}^{+\infty}\frac{\lambda_k^{1+2\beta+\frac{3}{2}+2\varepsilon}e^{-2t\lambda_k}}{1-e^{-2\lambda_k}} \vert(D_\varepsilon a)_k\vert^2\\ \nonumber
&=\frac{1}{t^{\frac{5}{2}+2\beta+2\varepsilon}}\sum_{k=1}^{+\infty}\frac{(t\lambda_k)^{(\frac{5}{2}+2\beta+2\varepsilon)}e^{-2t\lambda_k}} {1-e^{-2t\lambda_k}}\vert (D_\varepsilon a)_k\vert^2\\
 &  =\frac{1}{t^{\frac{5}{2}+2\beta+2\varepsilon}}\sum_{k=1}^{+\infty}\frac{(t\lambda_k)^{(\frac{5}{2}+2\beta+2\varepsilon)}} {e^{2t\lambda_k}-1}\vert (D_\varepsilon a)_k\vert^2\\
& \leq \frac{1}{t^{\frac{5}{2}+2\beta+2\varepsilon}} \sup_{x\geq 0} \frac{x^{\frac{5}{2}+2\beta+2\varepsilon}} {e^{x}-1}\sum_{k=1}^{+\infty}\vert (D_\varepsilon a)_k\vert^2\\
& \leq \frac{1}{t^{\frac{5}{2}+2\varepsilon}} \sup_{x\geq 0} \frac{x^{\frac{5}{2}+2\varepsilon}} {e^{x}-1}\vert (D_\varepsilon a)\vert^2 \leq C_0\frac{1}{t^{\frac{5}{2}+2\beta+2\varepsilon}}\vert  a\vert^2.
\end{align*}
So we can conclude that  as $t\rightarrow 0$, $\forall a\in \partial\,\calo$, $\vert Q_t^{-1/2}e^{tA_0}(-A_0D)a \vert^2$ blows up at most like $t^{-5/2-2\beta-2\varepsilon}$ and so the claim follows.
\qed

%
%
%

Now we introduce the operator $P$.
Let $\alpha>0$, let $v_1,..., v_n\in D((-A_0)^{\alpha})$ be linearly independent, and let $P$ be the projection on the span of $\<v_1,...,v_n\>$, namely
\begin{equation}\label{P-n-gen}
P:H\rightarrow H,\quad Px=\sum_{i=1}^n\<x,v_i\>v_i, \, \forall x \in H.
\end{equation}
We set moreover, noticing that $P=P^*$,
\begin{equation}\label{barQ_t-heat}
\bar Q_t:= PQ_tP=P (-A)^{-1-\beta}(I-e^{2tA})P
\end{equation}
Notice that
$P_\alpha:=(-A_0)^{\alpha}P$ is a continuous operator on $H$. Hence
$$
\overline{Pe^{tA_0}}B_0= P e^{tA_0}
(-A_0)^{\frac34+\eps}((-A_0)^{\frac14-\eps} D), \qquad
(\overline{Pe^{tA_0}}B_0)^*= ((-A_0)^{\frac14-\eps}D)^*
(-A_0)^{\frac34+\eps-\alpha} e^{tA_0}P_\alpha
$$
$$
\<Q_tP^*x,P^* x\>=
\<(I-e^{2tA_0}) (-A_0)^{-1-2\alpha-\beta} P_\alpha x,P_\alpha x\>
$$
The aim now is to verify that $\operatorname{Im}\overline{Pe^{tA_0}}(-A_0D)\subset \operatorname{Im}\bar Q_t^{1/2}$ and to estimate $\Vert \bar Q_t^{-1/2}Pe^{tA_0}(-A_0D)\Vert$.
\begin{lemma}\label{Lemma-barQt}
Let $\bar Q_t$ be defined in (\ref{barQ_t-heat}).
Let $\alpha>\beta + \frac14$. Then, for $\varepsilon\in(0,\frac14)$,
 \begin{equation}\label{barQ_t-norm-heat}
  \operatorname{Im}\overline{Pe^{tA}}(-A_0D)\subset \operatorname{Im}\bar Q_t^{1/2}, \quad \Vert \bar Q_t^{-1/2}\overline{Pe^{tA}}(-A_0D)\Vert \leq \frac{C}{t^{1-\varepsilon}}.
 \end{equation}
\end{lemma}

\dim In the proof we consider the case of the projection on the space generated by only one element $v\in D((-A)^\alpha)$, namely in the proof $P:H\rightarrow H,\; Px=\<x,v\>v,\, \forall x \in H$, the extension to a map as in \eqref{P-n-gen} being direct.
\newline We notice that
\[
\<\bar Q_t x,x\>=\vert\bar Q_t^{1/2}x\vert^2=\vert\<x,v\>\vert^2\sum_{k\geq1}\frac{1-e^{-2t\lambda_k}}{\lambda_k^{1+2\beta}}v_k^2
\]
and so
\begin{align}\label{stima-barQ_t}
\vert\bar Q_t^{-1/2}v\vert^2&=\frac{1}{\vert v\vert^2}\frac{1}{\sum_{k\geq1}\frac{1-e^{-2t\lambda_k}}{\lambda_k^{1+2\beta}}v_k^2}\leq
\frac{1}{t^{1+2\beta}\vert v\vert^2}\frac{1}{\sum_{k\geq1}\frac{e^{2t\lambda_k}-1}{(t\lambda_k)^{1+2\beta}e^{2t\lambda_k}}v_k^2}\\ \nonumber
&\leq \frac{1}{t^{1+2\beta}\vert v\vert^2}\frac{1}{\sum_{k\geq1}\frac{e^{2t\lambda_k}-1}{(t\lambda_k)^{1+2\beta}}v_k^2}
\leq \frac{1}{t^{1+2\beta}\vert v\vert^2}\frac{1}{\sum_{k\geq1}v_k^2}
\leq C_0  \frac{1}{t^{1+2\beta}\vert v\vert^4}
\end{align}
Moreover we can write, $\forall\, a \in \partial \,\calo$ and setting $D_\varepsilon :=(-A_0)^{\frac14-\eps}$
\[
\overline{Pe^{tA_0}}(-A_0D) a= Pe^{tA_0}(-A_0)^{\frac{3}{4}+\varepsilon}D_\varepsilon a
=\Big(\sum_{k\geq1}\lambda_k^{\frac{3}{4}+\varepsilon}e^{-t\lambda_k}(D_\varepsilon a)_kv_k\Big)\,v,
\]
so it is immediate to see that
$$ \operatorname{Im}\overline{Pe^{tA_0}}B\subset \operatorname{Im}\bar Q_t^{1/2}.$$
Taking into account that $D_\varepsilon$ is a bounded operator and that $v\in D((-A)^\alpha)$ we get
\begin{align*}
\Big(\sum_{k\geq1}\lambda_k^{\frac{3}{4}+\varepsilon}
e^{-t\lambda_k}(D_\varepsilon a)_kv_k\Big)^2&
= \Big(\sum_{k\geq1}\lambda_k^{\frac{3}{4}+\varepsilon-\alpha}e^{-t\lambda_k}(D_\varepsilon a)_k(\lambda_k)^\alpha v_k\Big)^2 \\
&=\vert \sum_{k\geq1}\lambda_k^{\frac{3}{4}+\varepsilon-\alpha}e^{-t\lambda_k}(D_\varepsilon a)_k(k^2)^\alpha v_k\vert^2 \\
&\leq
\left(\sum_{k\geq1}\lambda_k^{\frac{3}{2}+2\varepsilon-2\alpha}
e^{-2t\lambda_k}\vert (D_\varepsilon a)_k^2\right)
\left(\sum_{k\geq 1}\lambda_k^{2\alpha} v^2_k\right)
\\
&\leq  C_0 \frac{1}{t^{\frac{3}{2}+2\varepsilon-2\alpha}} \left(\sum_{k\geq1}(t\lambda_k)^{\frac{3}{2}+2\varepsilon-2\alpha}e^{-2t\lambda_k}\vert (D_\varepsilon a)_k\vert^2\right) \vert A^{\alpha} v\vert^2.
\end{align*}
So by these calculations, and by estimate \eqref{stima-barQ_t} we get, with the constant $C_0$ independent on $t$ and that may change value from line to line,
\begin{align*}
\vert \bar Q_t^{-1/2} \overline{Pe^{tA}}(-A_0D) a\vert^2
&\leq C_0 \frac{1}{t^{1+2\beta}\vert v\vert^4} \frac{1}{t^{\frac{3}{2}+2\varepsilon+2\beta-2\alpha}} \sum_{k\geq1}(t\lambda_k)^{\frac{1}{2}+2\varepsilon+2\beta-2\alpha}e^{-2t\lambda_k}\vert (D_\varepsilon a)_k^2\vert \vert A^{\alpha} v\vert^2\\
&\leq  C_0 \frac{1}{\vert v\vert^4} \frac{1}{t^{\frac{5}{2}+2\beta+2\varepsilon-2\alpha}} \sup_{x\geq 0}x^{\frac{3}{2}+2\varepsilon-2\alpha}e^{-2x}\vert D_\varepsilon a\vert^2.
\end{align*}
Choosing $\alpha>\frac{1}{4}+\beta=1-2\varepsilon$ we conclude the proof.
\qed

\subsection{The case of delay in the control}\label{sec-partialSmoothing-delay}

We consider equation \eqref{eq-astr} in Section \ref{SSE:DELAYEQUATION} and we notice that for every $(x_0,x_1) \in H$ the  covariance operator $Q_t$ of the stochastic convolution can be written as
\begin{equation}\label{eq:Q^0Q}
Q_t\left(
x_0, x_1
\right) =\left(
Q^0_t x_0, 0
 \right),
\end{equation}
where $Q^0_t$ is the selfadjoint operator in $\R^n$ defined as
\begin{equation}\label{eq:Q^0def}
Q^0_t:=\int_0^t e^{sa_0}\sigma\sigma^*
e^{sa_0^*}\, ds,
\end{equation}
see \cite{FGFM-I}, Lemma 4.6. So
$$
\operatorname{Im} Q_t
=\operatorname{Im} Q_t^0\times \left\lbrace 0\right\rbrace\subseteq \R^n
\times \left\lbrace 0\right\rbrace.$$
\begin{lemma} The operator $Q_t^0$ defined in (\ref{eq:Q^0def})
is invertible for all $t>0$ if and only if
\begin{equation}\label{cond-delay}
\operatorname{Im} (\sigma,a_0\sigma, \dots , a_0^{n-1}\sigma)=  \R^n,
\end{equation}
Let $P$ be the projection on the first component: $\forall (x_0,x_1)\in\, H$, $P(x_0,x_1)=(x_0,0)$.
{Then \eqref{eq:inclusionsmoothingC} holds if and only if
\begin{equation}\label{eq:inclusionsmoothingBbis}
\operatorname{Im}\left(e^{ta_0}b_0 +\int_{-d}^0 1_{[-t,0]}e^{(t+r)a_0}b_1(dr)
\right)\subseteq \operatorname{Im}(\sigma,a_0\sigma, \dots a_0^{n-1}\sigma ).
\end{equation}
If moreover
\begin{equation}\label{eq:hpdebregbis}
\operatorname{Im}\left(e^{ta_0}b_0 +\int_{-d}^0 1_{[-t,0]}e^{(t+r)a_0}b_1(dr)
\right)
\subseteq\operatorname{Im}\sigma,
\quad \forall t>0.
\end{equation}
then
\begin{equation}\label{eq:last-est}
\Vert (P Q_t P^*)^{-1/2}  \overline{Pe^{tA_1}}B_1\Vert \leq C
t^{-{1/2}}.
\end{equation}}
\end{lemma}
\dim
In this case \eqref{eq:inclusionsmoothingC} is written as
$$
\operatorname{Im}\overline{Pe^{tA_1}}B_1\subseteq
\operatorname{Im}(P Q_t P^*)^{1/2},\qquad \forall t>0.
$$
Recalling \eqref{eq:etAB}
we see that $\overline{Pe^{tA_1}}B_1$ is a bounded operator and can be written as
$$
\overline{Pe^{tA}}B_1=\left(
e^{ta_0 }x_0+\int_{-d}^{0}1_{[-t,0]} e^{(t+s)a_0 } b_1(ds),
0 \right):
$$
Hence the inclusion follows exactly as in \cite{FGFM-I} in the case when $B_1$ is a bounded operator: indeed, in \cite{FGFM-I}, formula (4.35), Proposition 4.11 the only difference is that in the present paper  $b_1(\cdot)$ is not necessarily absolutely continuous with respect to the Lebesgue measure) but this does not affect the image of
$\overline{Pe^{tA}}B_1$.
Hence, by  \cite{FGFM-I}, Proposition 4.11, it immediately follows \eqref{eq:last-est}.
\qed

\end{appendix}

%
%
%
%
%

\end{document}